\documentclass[11pt]{amsart}
\usepackage{enumitem}
\usepackage[paperheight=11in,paperwidth=8.5in,left=1in,top=1.25in,right=1in,bottom=1.25in,]{geometry}

\usepackage[super]{nth}
\usepackage[none]{hyphenat}
\usepackage[capitalize,nameinlink,noabbrev,nosort]{cleveref}
\usepackage{graphicx, amssymb}
\usepackage{tikz}   
\usetikzlibrary{arrows,calc,intersections,matrix,positioning,through}
\usepackage{tikz-cd}
\tikzset{commutative diagrams/.cd,every label/.append style = {font = \normalsize}}

\usepgflibrary{shapes.geometric}
\usetikzlibrary{shapes.geometric}
\DeclareMathOperator{\tree}{tree}
\DeclareMathOperator{\tab}{tab}
\DeclareMathOperator{\wt}{wt}
\DeclareMathOperator{\rt}{root}
\newcommand{\Sym}{\mathfrak{S}}
\newcommand{\Z}{\mathbb{Z}}

\newcommand{\swingleft}{\mathbin{{\rotatebox{180}{$\curvearrowright$}}}}

\newtheorem{theorem}{Theorem}[section]

\theoremstyle{definition}
\newtheorem{definition}[theorem]{Definition}
\newtheorem{example}[theorem]{Example}

\theoremstyle{remark}

\setlist[itemize]{leftmargin=*}
\setlist[enumerate]{leftmargin=*}
\newcommand{\ttt}{\tau}

\begin{document}
\begin{abstract}
The \emph{asymmetric simple exclusion process} (ASEP) 
is a model for translation in protein synthesis and traffic flow;
it can be defined as a Markov chain
describing particles hopping on a one-dimensional lattice.
In this article
I give an overview of some of the connections of 
the stationary distribution of the ASEP 
 to combinatorics (tableaux and multiline queues)
and special functions (Askey-Wilson polynomials,
Macdonald polynomials, and Schubert polynomials).  I also make some general 
observations about  positivity in Markov chains.
\end{abstract}

	\title[The combinatorics of hopping particles and positivity in Markov chains]{The combinatorics of hopping particles and positivity in Markov chains}
	\date{\today}
\author{Lauren K. Williams}
	\address{Department of Mathematics, 1 Oxford Street, Cambridge,
	MA 02138}
	\email{williams@math.harvard.edu}
\maketitle

\section{Introduction}

The goal of this article is to illustrate some of the elegant connections
between combinatorics and probability that arise when one studies Markov chains.
We will focus in particular on several variations of the 
\emph{asymmetric simple exclusion process}, illustrating combinatorial 
formulas for its stationary distribution and connections to special functions.
The last section of this article  makes some general observations about 
positivity in Markov chains, in the context of the 
 Markov Chain Tree Theorem.


The asymmetric simple exclusion process (or ASEP)
is a model for
 particles hopping on a one-dimensional lattice 
	(e.g. a \emph{line} or a \emph{ring}) such that
each site contains at most one particle.
The ASEP was introduced independently
in biology 
by Macdonald--Gibbs--Pipkin \cite{bio} 
and in mathematics by Spitzer \cite{Spitzer} around 1970, see also
\cite{Liggett}.
It exhibits boundary-induced phase transitions,
and has been cited as a model for translation in 
	protein synthesis,
sequence alignment, the nuclear pore complex, and traffic flow.

The ASEP has remarkable connections to a number of topics,
including the XXZ model \cite{Sandow},  vertex models  \cite{BorodinPetrov, BorodinWheeler},
	the Tracy-Widom distribution \cite{Johansson, TracyWidom}, 
	and the KPZ equation \cite{BG, CorwinShenTsai, CorwinShen, CorwinKnizel}.
	The ASEP is often viewed as a prototypical example of a random growth model from the 
	so-called KPZ universality class in $(1+1)$-dimensions, see 
	\cite{KPZ, Corwin, Quastel}.  
	However,
in this article we restrict our attention to the ASEP's relations
	to combinatorics 
	(including 
	staircase tableaux and multiline queues), as well as 
	to special functions
(including Askey-Wilson polynomials,
Macdonald polynomials,
and Schubert polynomials).   
Much of this article is based on 
joint works with Sylvie Corteel \cite{CW1, CW2,CWPNAS,CW4}, 
as well as Olya Mandelshtam \cite{CMW2}, 
and Donghun Kim \cite{KW}.


\section{The open boundary ASEP}

In the ASEP with open boundaries (see Figure 1), 
we have a one-dimensional
lattice of $n$ sites such that each site is either empty or occupied by 
a particle.  At most one particle may occupy a
given site.  
We can describe it informally as follows.
During each infinitesimal time interval $\mathit{dt}$, each particle at a site
$1\leq i \leq n-1$ 
has a probability $\mathit{dt}$ of jumping to the next
site on its right,
	provided it is empty, and each particle at a site $2\leq i \leq n$
	has 
a probability $q \mathit{dt}$ of jumping to the next site on its left,
	provided
	it is empty.  Furthermore, a particle
is added at site $i=1$ with probability $\alpha \mathit{dt}$ if site $1$
is empty and a particle is removed from site $n$ with probability
$\beta \mathit{dt}$ if this site is occupied.

More formally, we define this model as a 
discrete-time Markov chain.
\begin{figure}[h]
\resizebox{1.8in}{!}{
\begin{tikzpicture}
        \draw[black,thick](-.1,-.2)--(.1,-.2);
        \draw[black,thick](.25,-.2)--(.45,-.2);
        \draw[black,thick](.6,-.2)--(.8,-.2);
        \draw[black,thick](.95,-.2)--(1.15,-.2);
        \draw[black,thick](1.3,-.2)--(1.5,-.2);
        \draw[black,thick](1.65,-.2)--(1.85,-.2);
        \draw[black,thick](2,-.2)--(2.2,-.2);
        \filldraw[color=black,fill=white] (0,0) circle (3pt);
        \filldraw[color=black,fill=black] (.35,0) circle (3pt);
        \filldraw[color=black,fill=white] (.7,0) circle (3pt);
        \filldraw[color=black,fill=black] (1.05,0) circle (3pt);
        \filldraw[color=black,fill=black] (1.4,0) circle (3pt);
        \filldraw[color=black,fill=white] (1.75,0) circle (3pt);
        \filldraw[color=black,fill=black] (2.1,0) circle (3pt);
	 \draw (-.15,.4) node {\small{$\alpha$}};
	 \draw (-.15,.2) node {$\curvearrowright$};
	  \draw(.85,.2) node {$\curvearrowleft$};
	  \draw(.9,.42) node {\small{$q$}};
	  \draw(1.6,.2) node {$\curvearrowright$};
	  \draw(1.6,.42) node {\small{$1$}};
	  \draw(2.25,.2) node {$\curvearrowright$};
	  \draw(2.3,.42) node {\small{$\beta$}};
\end{tikzpicture}
}
	\caption{The (three-parameter) open boundary ASEP.}\label{fig:1} 
\end{figure}
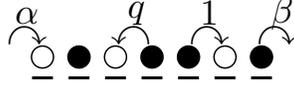

	\begin{definition}
	Let $\alpha$, $\beta$, and $q$ be constants  between $0$ and $1$.
Let $B_n$ be the set of all $2^n$ words of length $n$ in 
 $\{\circ, \bullet\}$.
The open boundary ASEP is the Markov chain on $B_n$ with
transition probabilities:
\begin{itemize}
	\item  If $\tau = A{\bullet}{\circ} B$ and
		$\sigma = A {\circ}{ \bullet} B$ (where $A$ and $ B$ are words
		in $\{\circ,\bullet\}$), then we have that 
		$\Pr(\tau\to \sigma) = \frac{1}{n+1}$ 
		and
		$\Pr(\sigma\to \tau) = \frac{q}{n+1}$ (particle hops right or left).
\item  If $\tau = \circ B$ and $\sigma = \bullet B$
	then $\Pr(\tau\to \sigma) = \frac{\alpha}{n+1}$ (particle enters the lattice from left).
\item  If $\tau = B \bullet$ and $\sigma = B \circ$
	then $\Pr(\tau\to \sigma) = \frac{\beta}{n+1}$ (particle exits the lattice to the right).
\item  Otherwise $\Pr(\tau\to \sigma) = 0$ for $\sigma \neq \tau$
	and $\Pr(\tau\to \tau) = 1 - \sum_{\sigma \neq \tau} \Pr(\tau\to \sigma)$.
\end{itemize}
\end{definition}


In the long time limit, the system reaches a steady state where all
the probabilities ${\pi}(\ttt)$ of finding
the system in configuration $\ttt$ are
stationary, i.e.\ satisfy
$	\frac{d}{dt} {\pi}(\tau) = 0.$
Moreover, the stationary distribution is unique.  We can compute it 
 by solving for the left eigenvector of the transition
matrix with eigenvalue $1$, or equivalently, by solving the 
\emph{global balance} equations: for all states $\tau\in B_n$, we have
\begin{equation*}
	\pi(\tau) \sum_{\sigma\neq \tau} \Pr(\tau\to \sigma) = 
	\sum_{\sigma\neq \tau} \pi(\sigma) \Pr(\sigma\to \tau),
\end{equation*}
where both sums are over all states $\sigma\neq \tau$.

The steady state probabilities are rational expressions in 
$\alpha, \beta$ and $q$.  For convenience, we clear denominators,
obtaining ``unnormalized probabilities'' $\Psi(\tau)$ 
 which are equal to the $\pi(\tau)$ up to a constant: that is, 
$\pi(\tau)=\frac{\Psi(\tau)}{Z_n}$, where $Z_n=Z_n(\alpha, \beta, q)$ 
is the \emph{partition function} $\sum_{\tau\in B_n} \Psi(\tau)$.
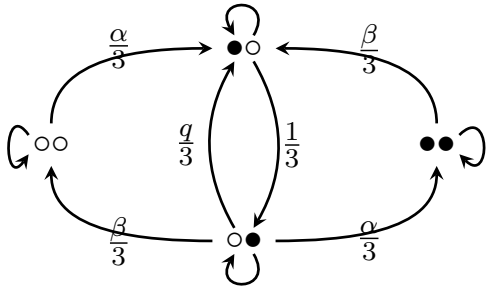
\begin{figure}[h]
\resizebox{2.7in}{!}{
\begin{tikzpicture}
	\draw (0,0) node (bw) {$\bullet \circ$};
	\draw (-2,-1) node (ww) {$\circ \circ$};
	\draw (-2.1,-.9) node (ww+) {};
	\draw (-2.1,-1.1) node (ww-) {};
	\draw (2,-1) node (bb) {$\bullet \bullet$};
	\draw (2.1,-.9) node (bb+) {};
	\draw (2.1,-1.1) node (bb-) {};
	\draw (0,-2) node (wb) {$\circ \bullet$};
	\draw (.1,-2) node (wb+) {};
	\draw (-.1,-2) node (wb-) {};
	\draw(.1,-.0) node (bw+) {};
	\draw(-.1,-.0) node (bw-) {};
        \draw (1.3, 0) node {$\frac{\beta}{3}$};
        \draw (-1.3, 0) node {$\frac{\alpha}{3}$};
        \draw (1.3, -2) node {$\frac{\alpha}{3}$};
        \draw (-1.3, -2) node {$\frac{\beta}{3}$};
        \draw (.5, -1) node {$\frac{1}{3}$};
        \draw (-.6, -1) node {$\frac{q}{3}$};

	\draw[thick, >=stealth, ->] (ww.north) .. controls +(up:7mm) and +(left:7mm) .. (bw.west);
	\draw[thick, >=stealth, ->] (bw+.south) to [out=300, in=60] (wb+.north);
	\draw[thick, >=stealth, ->] (wb-.north) to [out=120, in=240] (bw-.south);
	\draw[thick, >=stealth, ->] (bw+.north) to [out=60, in=120, looseness=6] (bw-.north);
	\draw[thick, >=stealth, ->] (wb+.south) to [out=-60, in=-120, looseness=6] (wb-.south);
	\draw[thick, >=stealth, ->] (bb+.east) to [out=45, in=-45, looseness=6] (bb-.east);
	\draw[thick, >=stealth, ->] (ww+.west) to [out=125, in=-125, looseness=6] (ww-.west);
	\draw[thick, >=stealth, ->] (bb.north) .. controls +(up:7mm) and +(right:7mm) .. (bw.east);
	\draw[thick, >=stealth, ->] (wb.west) .. controls +(left:7mm) and +(down:7mm) .. (ww.south);
	\draw[thick, >=stealth, ->] (wb.east) .. controls +(right:7mm) and +(down:7mm) .. (bb.south);
\end{tikzpicture}
}
	\caption{The state diagram of the open-boundary 
	ASEP on a lattice of $2$ sites.}\label{fig:2}
\end{figure}

\begin{example} \label{ex:1}
	\cref{fig:2}
	shows the state diagram of the open-boundary 
ASEP when $n=2$, 
	and \cref{table:1} gives the 
corresponding unnormalized probabilities.
Therefore we have 
$\pi(\bullet \bullet) = \frac{\alpha^2}{Z_2}$,
$\pi(\bullet \circ) = \frac{\alpha \beta(\alpha+\beta+q)}{Z_2}$, 
$\pi(\circ \bullet) = \frac{\alpha \beta}{Z_2}$, and 
$\pi(\circ \circ) = \frac{\beta^2}{Z_2}$, where
$Z_2 = \alpha^2+\alpha \beta(\alpha+\beta+q) + \alpha \beta + \beta^2$.
\end{example}
\begin{table}[h]
\begin{center}
\begin{tabular}{|c c| }
    \hline
	State $\tau$ & Unnormalized probability $\Psi(\tau)$\\
    \hline 
    $\bullet \bullet$ & $\alpha^2$\\
	$\bullet \circ$  & $\alpha \beta(\alpha+\beta+q)$\\
	$\circ \bullet$  & $\alpha \beta$\\
	$\circ \circ$ & $\beta^2$\\
    \hline
    \end{tabular}
\end{center}
	\caption{Probabilities for the open-boundary
	ASEP on a lattice of $2$ sites.}\label{table:1}
\end{table}

For $n=3$, 
if we again write 
each probability $\pi(\tau)=\frac{\Psi(\tau)}{Z_3}$,
we find that 
$Z_3(\alpha,\beta,q)$ is a polynomial which is 
\emph{manifestly positive} -- that it, it has 
 only positive coefficients.  Also, $Z_3$ has
$24$ terms (counted with multiplicity):  
$Z_3(1, 1, 1) = 24$.
Computing more examples quickly leads to the conjecture
that the {partition function}
$Z_n=Z_n(\alpha,\beta,q)$ is a (manifestly) positive polynomial with 
$(n-1)$! terms.

In algebraic combinatorics,  if a quantity of interest is known or believed
to be a positive integer or a polynomial with positive coefficients,  one seeks an 
interpretation of this quantity as counting some combinatorial objects.  For example,
one seeks to express such a polynomial as a generating function for certain
tableaux or graphs or permutations, etc. 
A prototypical example is the Schur polynomial $s_{\lambda}(x_1,\dots,x_n)$ \cite{EC2}: there are several
formulas for it, including the bialternant formula and the Jacobi-Trudi formula, but both of 
these involve determinants and neither 
 makes it obvious that the Schur polynomial has positive coefficients.  However, one can express 
the Schur polynomial as the generating function for semistandard tableaux of shape $\lambda$,
and this formula makes manifest the positivity of coefficients \cite{EC2}.

Given the above, and 
our observations on the positivity of the partition function $Z_n(\alpha,\beta,q)$,
the natural question is:
can we express each probability
as a (manifestly positive) sum over some set of combinatorial objects?
We will explain how to answer this question using (a special case of) 
the staircase
tableaux of \cite{CW4}.

\subsection{The open-boundary ASEP and $\alpha \beta$-staircase tableaux}
In what follows,
we will depict Young diagrams
in Russian notation (with the corner at the bottom).  

\begin{definition}
		An  \emph{$\alpha \beta$-staircase tableau} $T$ of size $n$ is a 
Young diagram of shape $(n, n-1, \dots, 2, 1)$ 
(drawn in Russian notation)
such that each box is either empty
or contains an $\alpha$ or $\beta$, such that:
\begin{enumerate}
\item no box in the top row is empty
\item each box southeast of a $\beta$ and in the same diagonal
 as that $\beta$ is empty.
 \item each box southwest of an $\alpha$ and in the same diagonal
 as that $\alpha$ is empty.
\end{enumerate}
\end{definition}

See Figure 3.
It is an exercise to verify that 
there are  $(n+1)!$ $\alpha \beta$-staircase tableaux
of size $n$.

	\begin{figure}[h]
\resizebox{4.2in}{!}{
\begin{tikzpicture}
        \draw[black,thick](0,.2)--(1.4,-1.2);
        \draw[black,thick](.4,.2)--(1.6,-1);
        \draw[black,thick](.8,.2)--(1.8,-.8);
        \draw[black,thick](1.2,.2)--(2,-.6);
        \draw[black,thick](1.6,.2)--(2.2,-.4);
        \draw[black,thick](2,.2)--(2.4,-.2);
        \draw[black,thick](2.4,.2)--(2.6,0);
        \draw[black,thick](0,.2)--(-.2,0);
        \draw[black,thick](.4,.2)--(0,-.2);
        \draw[black,thick](.8,.2)--(.2,-.4);
        \draw[black,thick](1.2,.2)--(.4,-.6);
        \draw[black,thick](1.6,.2)--(.6,-.8);
        \draw[black,thick](2,.2)--(.8,-1);
        \draw[black,thick](2.4,.2)--(1,-1.2);
        \draw[black,thick](2.6,0)--(1.2,-1.4);
        \draw[black,thick](-.2,0)--(1.2,-1.4);
	\draw (0,.4) node {$\circ$};
	\draw (.4,.4) node {$\bullet$};
	\draw (.8,.4) node {$\circ$};
	\draw (1.2,.4) node {$\bullet$};
	\draw (1.6,.4) node {$\bullet$};
	\draw (2,.4) node {$\circ$};
	\draw (2.4,.4) node {$\bullet$};
	 \draw (0,0) node {\scriptsize{$\beta$}};
	 \draw (.4,0) node {\scriptsize{$\alpha$}};
	 \draw (.8,0) node {\scriptsize{$\beta$}};
	 \draw (1.2,0) node {\scriptsize{$\alpha$}};
	 \draw (1.2,-.8) node {\scriptsize{$\alpha$}};
	 \draw (1.6,0) node {\scriptsize{$\alpha$}};
	 \draw (1.8,-.2) node {\scriptsize{$\beta$}};
	 \draw (2,0) node {\scriptsize{$\beta$}};
	 \draw (2.4,0) node {\scriptsize{$\alpha$}};
\end{tikzpicture}\hspace{1in}
\begin{tikzpicture}
        \draw[black,thick](0,.2)--(1.4,-1.2);
        \draw[black,thick](.4,.2)--(1.6,-1);
        \draw[black,thick](.8,.2)--(1.8,-.8);
        \draw[black,thick](1.2,.2)--(2,-.6);
        \draw[black,thick](1.6,.2)--(2.2,-.4);
        \draw[black,thick](2,.2)--(2.4,-.2);
        \draw[black,thick](2.4,.2)--(2.6,0);
        \draw[black,thick](0,.2)--(-.2,0);
        \draw[black,thick](.4,.2)--(0,-.2);
        \draw[black,thick](.8,.2)--(.2,-.4);
        \draw[black,thick](1.2,.2)--(.4,-.6);
        \draw[black,thick](1.6,.2)--(.6,-.8);
        \draw[black,thick](2,.2)--(.8,-1);
        \draw[black,thick](2.4,.2)--(1,-1.2);
        \draw[black,thick](2.6,0)--(1.2,-1.4);
        \draw[black,thick](-.2,0)--(1.2,-1.4);
	 \draw (0,0) node {\scriptsize{$\beta$}};
	 \draw (.4,0) node {\scriptsize{$\alpha$}};
	 \draw (.6,-.2) node {\scriptsize{$q$}};
	 \draw (.8,0) node {\scriptsize{$\beta$}};
	 \draw (1.2,0) node {\scriptsize{$\alpha$}};
	 \draw (1.2,-.8) node {\scriptsize{$\alpha$}};
	 \draw (1.6,0) node {\scriptsize{$\alpha$}};
	 \draw (1.6,-.4) node {\scriptsize{$q$}};
	 \draw (1.8,-.2) node {\scriptsize{$\beta$}};
	 \draw (2,0) node {\scriptsize{$\beta$}};
	 \draw (2.4,0) node {\scriptsize{$\alpha$}};
\end{tikzpicture}
}
\caption{At left: an $\alpha\beta$-staircase tableau $T$
		of type $(\circ \bullet \circ \bullet \bullet
		\circ \bullet)$. At right: $T$
		with a $q$ in each unrestricted box. 
		We have $\wt(T)=\alpha^5 \beta^4 q^2$.} \label{fig:3}
\end{figure}
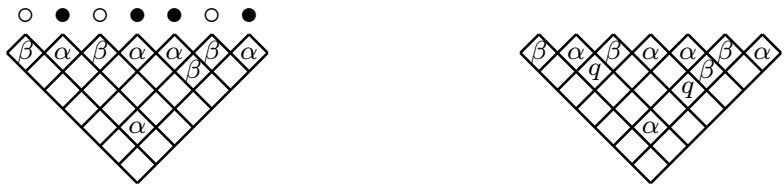

\begin{definition}
Some boxes in a tableau
	are forced to be empty because of conditions (2) or (3) above;
 we refer
to all other empty boxes as \emph{unrestricted}.  (The
unrestricted
boxes are  those whose nearest neighbors on the diagonals
	to the northwest and northeast, respectively, are an $\alpha$ and $\beta$.)

	After placing a $q$
	in each unrestricted box, we define the 
	\emph{weight} $\wt(T)$ of $T$ to be 
	$\alpha^i \beta^j q^k$ where $i, j$ and $k$ are the 
	numbers of $\alpha$'s, $\beta$'s, and $q$'s in $T$.

	The \emph{type} of $T$ is the word
	obtained by reading the letters in the top row of $T$
	and replacing each $\alpha$ by $\bullet$ and $\beta$ by 
	$\circ$, 
see Figure 3.
\end{definition}

The following result \cite{CW1, CW2, CW4}
gives a combinatorial formula for the steady state 
probabilities of the ASEP.  (Note that \cite{CW1, CW2} used 
\emph{permutation tableaux}, which are in bijection with staircase tableaux,
and are closely connected to the \emph{positive Grassmannian}.)
The $q=0$ case was previously
studied by Duchi--Schaeffer in \cite{DS}.

\begin{theorem}\label{thm:1}
Consider the ASEP with open boundaries on a lattice of $n$ sites.
Let $\tau=(\tau_1,\dots, \tau_n)\in \{\bullet,\circ\}^n$ be
	a state.  Then the 
	unnormalized steady state probability $\Psi(\tau)$ is equal to 
 $\sum_T \wt(T)$, where the sum is over the 
$\alpha \beta$-staircase tableaux of type $\tau$.

Equivalently, if 
we let $\mathcal{T}_n$ be the set of $\alpha \beta$-staircase 
tableaux of size $n$, and  $Z_n:= \sum_{T\in \mathcal{T}_n} \wt(T)$
be the weight generating function for these tableaux, then 
the	steady state probability $\pi(\tau)$
is $\frac{\sum_T \wt(T)}{Z_n}$, where the sum is over the 
$\alpha \beta$-staircase tableaux of type $\tau$.
\end{theorem}

In the case $n=2$, there are six tableaux of size
$2$, shown in Figure 4 and arranged by type.   Computing the weights
of the tableaux of the various types reproduces the 
results from \cref{ex:1}.

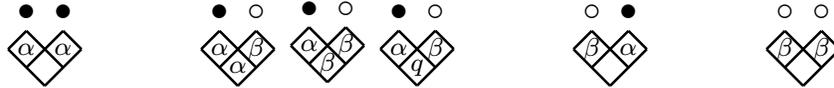
\begin{figure}[h]
\resizebox{4.5in}{!}{
\begin{tikzpicture}
        \draw[black,thick](0,.2)--(.4,-.2);
        \draw[black,thick](.4,.2)--(.6,0);
        \draw[black,thick](0,.2)--(-.2,0);
        \draw[black,thick](.4,.2)--(0,-.2);
        \draw[black,thick](.6,0)--(.2,-.4);
        \draw[black,thick](-.2,0)--(0.2,-0.4);
	\draw (0,.4) node {$\bullet$};
	\draw (.4,.4) node {$\bullet$};
	 \draw (0,0) node {\scriptsize{$\alpha$}};
	 \draw (.4,0) node {\scriptsize{$\alpha$}};
\end{tikzpicture}
	\hspace{1cm}
\begin{tikzpicture}
        \draw[black,thick](0,.2)--(.4,-.2);
        \draw[black,thick](.4,.2)--(.6,0);
        \draw[black,thick](0,.2)--(-.2,0);
        \draw[black,thick](.4,.2)--(0,-.2);
        \draw[black,thick](.6,0)--(.2,-.4);
        \draw[black,thick](-.2,0)--(0.2,-0.4);
	\draw (0,.4) node {$\bullet$};
	\draw (.4,.4) node {$\circ$};
	 \draw (0,0) node {\scriptsize{$\alpha$}};
	 \draw (.2,-.2) node {\scriptsize{$\alpha$}};
	 \draw (.4,0) node {\scriptsize{$\beta$}};
\end{tikzpicture}
\begin{tikzpicture}
        \draw[black,thick](0,.2)--(.4,-.2);
        \draw[black,thick](.4,.2)--(.6,0);
        \draw[black,thick](0,.2)--(-.2,0);
        \draw[black,thick](.4,.2)--(0,-.2);
        \draw[black,thick](.6,0)--(.2,-.4);
        \draw[black,thick](-.2,0)--(0.2,-0.4);
	\draw (0,.4) node {$\bullet$};
	\draw (.4,.4) node {$\circ$};
	 \draw (0,0) node {\scriptsize{$\alpha$}};
	 \draw (.2,-.2) node {\scriptsize{$\beta$}};
	 \draw (.4,0) node {\scriptsize{$\beta$}};
\end{tikzpicture}
\begin{tikzpicture}
        \draw[black,thick](0,.2)--(.4,-.2);
        \draw[black,thick](.4,.2)--(.6,0);
        \draw[black,thick](0,.2)--(-.2,0);
        \draw[black,thick](.4,.2)--(0,-.2);
        \draw[black,thick](.6,0)--(.2,-.4);
        \draw[black,thick](-.2,0)--(0.2,-0.4);
	\draw (0,.4) node {$\bullet$};
	\draw (.4,.4) node {$\circ$};
	 \draw (0,0) node {\scriptsize{$\alpha$}};
	 \draw (.2,-.2) node {\scriptsize{$q$}};
	 \draw (.4,0) node {\scriptsize{$\beta$}};
\end{tikzpicture}
	\hspace{1cm}
\begin{tikzpicture}
        \draw[black,thick](0,.2)--(.4,-.2);
        \draw[black,thick](.4,.2)--(.6,0);
        \draw[black,thick](0,.2)--(-.2,0);
        \draw[black,thick](.4,.2)--(0,-.2);
        \draw[black,thick](.6,0)--(.2,-.4);
        \draw[black,thick](-.2,0)--(0.2,-0.4);
	\draw (0,.4) node {$\circ$};
	\draw (.4,.4) node {$\bullet$};
	 \draw (0,0) node {\scriptsize{$\beta$}};
	 \draw (.4,0) node {\scriptsize{$\alpha$}};
\end{tikzpicture}
	\hspace{1cm}
\begin{tikzpicture}
        \draw[black,thick](0,.2)--(.4,-.2);
        \draw[black,thick](.4,.2)--(.6,0);
        \draw[black,thick](0,.2)--(-.2,0);
        \draw[black,thick](.4,.2)--(0,-.2);
        \draw[black,thick](.6,0)--(.2,-.4);
        \draw[black,thick](-.2,0)--(0.2,-0.4);
	\draw (0,.4) node {$\circ$};
	\draw (.4,.4) node {$\circ$};
	 \draw (0,0) node {\scriptsize{$\beta$}};
	 \draw (.4,0) node {\scriptsize{$\beta$}};
\end{tikzpicture}
}
\caption{The six $\alpha\beta$-staircase tableau $T$
	of size $2$.} 
	\label{fig:4}
\end{figure}

\subsection{Lumping of Markov chains}
Given  a combinatorial formula such as \cref{thm:1}, how can we prove it? 
One option is to realize the ASEP as a \emph{lumping} (or \emph{projection}) of 
a Markov chain on tableaux \cite{CW2}.
(See also \cite{DS} for the case $q=0$.)
Recall that 
we have a surjection
$f:\mathcal{T}_n \to B_n$, which maps an 
 $\alpha \beta$-staircase tableau
to its \emph{type}.
We'd like to construct a Markov chain
on tableaux whose projection via $f$ recovers the ASEP.
If we can do so, and moreover show that 
	  the steady state probability $\pi(T)$
		is proportional to $\wt(T)$, then 
		we will have proved \cref{thm:1}.



	\begin{definition}
Let $\{X_t\}$ be a Markov chain on state space $\Omega_X$ with transition matrix $P$,
and let $f:\Omega_X \to \Omega_Y$ be a surjective map.  
Suppose there is an $|\Omega_Y| \times |\Omega_Y|$ matrix $Q$ such that 
	for all $y_0, y_1\in \Omega_Y$, if $f(x_0)=y_0$, then 
	\begin{equation} \label{eq:lumping}
	\sum_{x: f(x)=y_1} P(x_0,x) = Q(y_0, y_1).
	\end{equation}
Then $\{f(X_t)\}$  is a Markov chain on 
$\Omega_Y$ with transition matrix $Q$.
We say that $\{f(X_t)\}$ is a 
	\emph{(strong) lumping} 
	of $\{X_t\}$ and 
$\{X_t\}$ is a 
	\emph{(strong) lift} of $\{f(X_t)\}$.

Suppose	$\pi$ is a stationary distribution
for $\{X_t\}$, and 
	let $\pi_f$ be the measure on $\Omega_Y$ defined by 
	$\pi_f(y) = \sum_{x: f(x)=y} \pi(x)$. Then $\pi_f$
	is a stationary distribution for $\{f(X_t)\}$.
	\end{definition}



See \cite{KemenySnell, Pang} for a thorough discussion of 
lumping.

The ASEP can be lifted to a Markov chain on $\alpha \beta$-staircase
tableaux \cite{CW2}, see Figure 5.  
In each diagram, the grey boxes
represent boxes that must be empty.  Note
that the 
remaining empty boxes on the left and right side of a ``$\mapsto$''
are in bijection with each other; they must be filled the 
same way.
The lifted chain has the 
particularly nice property that the left hand side of \eqref{eq:lumping} 
always has at most one nonzero term.


	\begin{figure}
\resizebox{2.8in}{!}{
 \begin{tikzpicture}
	\draw (.2,-.2) node [draw,scale=.7,diamond,fill=gray!30]{};
	\draw (.4,-.4) node [draw,scale=.7,diamond,fill=gray!30]{};
	\draw (.6,-.6) node [draw,scale=.7,diamond,fill=gray!30]{};
	\draw (.8,-.8) node [draw,scale=.7,diamond,fill=gray!30]{};
	\draw (1,-1) node [draw,scale=.7,diamond,fill=gray!30]{};
	\draw (1.2,-1.2) node [draw,scale=.7,diamond,fill=gray!30]{};
        \draw[black,thick](0,.2)--(1.4,-1.2);
        \draw[black,thick](.4,.2)--(1.6,-1);
        \draw[black,thick](.8,.2)--(1.8,-.8);
        \draw[black,thick](1.2,.2)--(2,-.6);
        \draw[black,thick](1.6,.2)--(2.2,-.4);
        \draw[black,thick](2,.2)--(2.4,-.2);
        \draw[black,thick](2.4,.2)--(2.6,0);
        \draw[black,thick](0,.2)--(-.2,0);
        \draw[black,thick](.4,.2)--(0,-.2);
        \draw[black,thick](.8,.2)--(.2,-.4);
        \draw[black,thick](1.2,.2)--(.4,-.6);
        \draw[black,thick](1.6,.2)--(.6,-.8);
        \draw[black,thick](2,.2)--(.8,-1);
        \draw[black,thick](2.4,.2)--(1,-1.2);
        \draw[black,thick](2.6,0)--(1.2,-1.4);
        \draw[black,thick](-.2,0)--(1.2,-1.4);
	\draw (0,.4) node {$\circ$};
	\draw (.4,.4) node {$\bullet$};
	\draw (.8,.4) node {$\bullet$};
	\draw (1.2,.4) node {$\bullet$};
	\draw (1.6,.4) node {$\circ$};
	\draw(-.2,.6) node {{$\curvearrowright$}};
	 \draw (0,0) node {\scriptsize{$\beta$}};
	 \draw (.4,0) node {\scriptsize{$\alpha$}};
	 \draw (.8,0) node {\scriptsize{$\alpha$}};
	 \draw (1.2,0) node {\scriptsize{$\alpha$}};
	 \draw (1.6,0) node {\scriptsize{$\beta$}};
	 \draw (2.2,.4) node {\scriptsize{$\cdots$}};
	 \draw (2.8,-.7) node {$\alpha$};
	 \draw (2.8,-.9) node {$\mapsto$};
\end{tikzpicture}
 \begin{tikzpicture}
	\draw (.6,-.6) node [draw,scale=.7,diamond,fill=gray!30]{};
	\draw (.8,-.4) node [draw,scale=.7,diamond,fill=gray!30]{};
	\draw (1,-.2) node [draw,scale=.7,diamond,fill=gray!30]{};
	\draw (1.6,-.4) node [draw,scale=.7,diamond,fill=gray!30]{};
	\draw (1.8,-.6) node [draw,scale=.7,diamond,fill=gray!30]{};
        \draw[black,thick](0,.2)--(1.4,-1.2);
        \draw[black,thick](.4,.2)--(1.6,-1);
        \draw[black,thick](.8,.2)--(1.8,-.8);
        \draw[black,thick](1.2,.2)--(2,-.6);
        \draw[black,thick](1.6,.2)--(2.2,-.4);
        \draw[black,thick](2,.2)--(2.4,-.2);
        \draw[black,thick](2.4,.2)--(2.6,0);
        \draw[black,thick](0,.2)--(-.2,0);
        \draw[black,thick](.4,.2)--(0,-.2);
        \draw[black,thick](.8,.2)--(.2,-.4);
        \draw[black,thick](1.2,.2)--(.4,-.6);
        \draw[black,thick](1.6,.2)--(.6,-.8);
        \draw[black,thick](2,.2)--(.8,-1);
        \draw[black,thick](2.4,.2)--(1,-1.2);
        \draw[black,thick](2.6,0)--(1.2,-1.4);
        \draw[black,thick](-.2,0)--(1.2,-1.4);
	\draw (0,.4) node {$\bullet$};
	\draw (.4,.4) node {$\bullet$};
	\draw (.8,.4) node {$\bullet$};
	\draw (1.2,.4) node {$\bullet$};
	\draw (1.6,.4) node {$\circ$};
	 \draw (0,0) node {\scriptsize{$\alpha$}};
	 \draw (.4,0) node {\scriptsize{$\alpha$}};
	 \draw (.8,0) node {\scriptsize{$\alpha$}};
	 \draw (1.2,0) node {\scriptsize{$\alpha$}};
	 \draw (1.4,-.2) node {\scriptsize{$\beta$}};
	 \draw (1.6,0) node {\scriptsize{$\beta$}};
	 \draw (2.2,.4) node {\scriptsize{$\cdots$}};
\end{tikzpicture}
}
\hspace{1cm} \resizebox{2.8in}{!}{
 \begin{tikzpicture}
	\draw (.6,-.6) node [draw,scale=.7,diamond,fill=gray!30]{};
	\draw (.8,-.4) node [draw,scale=.7,diamond,fill=gray!30]{};
	\draw (1.4,-.2) node [draw,scale=.7,diamond,fill=gray!30]{};
	\draw (1.6,-.4) node [draw,scale=.7,diamond,fill=gray!30]{};
	\draw (1.8,-.6) node [draw,scale=.7,diamond,fill=gray!30]{};
        \draw[black,thick](0,.2)--(1.4,-1.2);
        \draw[black,thick](.4,.2)--(1.6,-1);
        \draw[black,thick](.8,.2)--(1.8,-.8);
        \draw[black,thick](1.2,.2)--(2,-.6);
        \draw[black,thick](1.6,.2)--(2.2,-.4);
        \draw[black,thick](2,.2)--(2.4,-.2);
        \draw[black,thick](2.4,.2)--(2.6,0);
        \draw[black,thick](0,.2)--(-.2,0);
        \draw[black,thick](.4,.2)--(0,-.2);
        \draw[black,thick](.8,.2)--(.2,-.4);
        \draw[black,thick](1.2,.2)--(.4,-.6);
        \draw[black,thick](1.6,.2)--(.6,-.8);
        \draw[black,thick](2,.2)--(.8,-1);
        \draw[black,thick](2.4,.2)--(1,-1.2);
        \draw[black,thick](2.6,0)--(1.2,-1.4);
        \draw[black,thick](-.2,0)--(1.2,-1.4);
	\draw (0,.4) node {$\bullet$};
	\draw (.4,.4) node {$\circ$};
	\draw(1,.6) node {{$\curvearrowright$}};
	\draw (.8,.4) node {$\bullet$};
	\draw (1.2,.4) node {$\circ$};
	\draw (1.8,.4) node {$\cdots$};
	 \draw (0,0) node {\scriptsize{$\alpha$}};
	 \draw (.4,0) node {\scriptsize{$\beta$}};
	 \draw (.8,0) node {\scriptsize{$\alpha$}};
	 \draw (1,-.2) node {\scriptsize{$\alpha$}};
	 \draw (1.2,0) node {\scriptsize{$\beta$}};
	 \draw (2.8,-.7) node {$1$};
	 \draw (2.8,-.9) node {$\mapsto$};
\end{tikzpicture}

	\begin{tikzpicture}
	\draw (.6,-.2) node [draw,scale=.7,diamond,fill=gray!30]{};
	\draw (.8,-.4) node [draw,scale=.7,diamond,fill=gray!30]{};
	\draw (1,-.6) node [draw,scale=.7,diamond,fill=gray!30]{};
	\draw (1.2,-.8) node [draw,scale=.7,diamond,fill=gray!30]{};
	\draw (1.4,-1) node [draw,scale=.7,diamond,fill=gray!30]{};
        \draw[black,thick](0,.2)--(1.4,-1.2);
        \draw[black,thick](.4,.2)--(1.6,-1);
        \draw[black,thick](.8,.2)--(1.8,-.8);
        \draw[black,thick](1.2,.2)--(2,-.6);
        \draw[black,thick](1.6,.2)--(2.2,-.4);
        \draw[black,thick](2,.2)--(2.4,-.2);
        \draw[black,thick](2.4,.2)--(2.6,0);
        \draw[black,thick](0,.2)--(-.2,0);
        \draw[black,thick](.4,.2)--(0,-.2);
        \draw[black,thick](.8,.2)--(.2,-.4);
        \draw[black,thick](1.2,.2)--(.4,-.6);
        \draw[black,thick](1.6,.2)--(.6,-.8);
        \draw[black,thick](2,.2)--(.8,-1);
        \draw[black,thick](2.4,.2)--(1,-1.2);
        \draw[black,thick](2.6,0)--(1.2,-1.4);
        \draw[black,thick](-.2,0)--(1.2,-1.4);
	\draw (0,.4) node {$\bullet$};
	\draw (.4,.4) node {$\circ$};
	\draw (.8,.4) node {$\circ$};
	\draw (1.2,.4) node {$\bullet$};
	\draw (1.8,.4) node {$\cdots$};
	 \draw (0,0) node {\scriptsize{$\alpha$}};
	 \draw (.2,-.2) node {\scriptsize{$\alpha$}};
	 \draw (.4,0) node {\scriptsize{$\beta$}};
	 \draw (.8,0) node {\scriptsize{$\beta$}};
	 \draw (1.2,0) node {\scriptsize{$\alpha$}};
\end{tikzpicture}
}

\vspace{.3cm}
\resizebox{2.8in}{!}{
		\begin{tikzpicture}
	\draw (1.2,-1.2) node [draw,scale=.7,diamond,fill=gray!30]{};
	\draw (1.4,-1) node [draw,scale=.7,diamond,fill=gray!30]{};
	\draw (1.6,-.8) node [draw,scale=.7,diamond,fill=gray!30]{};
	\draw (1.8,-.6) node [draw,scale=.7,diamond,fill=gray!30]{};
	\draw (2,-.4) node [draw,scale=.7,diamond,fill=gray!30]{};
	\draw (2.2,-.2) node [draw,scale=.7,diamond,fill=gray!30]{};
        \draw[black,thick](0,.2)--(1.4,-1.2);
        \draw[black,thick](.4,.2)--(1.6,-1);
        \draw[black,thick](.8,.2)--(1.8,-.8);
        \draw[black,thick](1.2,.2)--(2,-.6);
        \draw[black,thick](1.6,.2)--(2.2,-.4);
        \draw[black,thick](2,.2)--(2.4,-.2);
        \draw[black,thick](2.4,.2)--(2.6,0);
        \draw[black,thick](0,.2)--(-.2,0);
        \draw[black,thick](.4,.2)--(0,-.2);
        \draw[black,thick](.8,.2)--(.2,-.4);
        \draw[black,thick](1.2,.2)--(.4,-.6);
        \draw[black,thick](1.6,.2)--(.6,-.8);
        \draw[black,thick](2,.2)--(.8,-1);
        \draw[black,thick](2.4,.2)--(1,-1.2);
        \draw[black,thick](2.6,0)--(1.2,-1.4);
        \draw[black,thick](-.2,0)--(1.2,-1.4);
	\draw (.2,.4) node {$\cdots$};
	\draw (.8,.4) node {$\bullet$};
	\draw (1.2,.4) node {$\circ$};
	\draw (1.6,.4) node {$\circ$};
	\draw (2,.4) node {$\circ$};
	\draw (2.4,.4) node {$\bullet$};
	\draw(2.6,.6) node {{$\curvearrowright$}};
	 \draw (.8,0) node {\scriptsize{$\alpha$}};
	 \draw (1.2,0) node {\scriptsize{$\beta$}};
	 \draw (1.6,0) node {\scriptsize{$\beta$}};
	 \draw (2,0) node {\scriptsize{$\beta$}};
	 \draw (2.4,0) node {\scriptsize{$\alpha$}};
	\draw(-.2,.6) node {};
	 \draw (2.85,-.65) node {$\beta$};
	 \draw (2.8,-.9) node {$\mapsto$};
\end{tikzpicture}
 \begin{tikzpicture}
	\draw (.6,-.6) node [draw,scale=.7,diamond,fill=gray!30]{};
	\draw (.8,-.4) node [draw,scale=.7,diamond,fill=gray!30]{};
	\draw (1.4,-.2) node [draw,scale=.7,diamond,fill=gray!30]{};
	\draw (1.6,-.4) node [draw,scale=.7,diamond,fill=gray!30]{};
	\draw (1.8,-.6) node [draw,scale=.7,diamond,fill=gray!30]{};
        \draw[black,thick](0,.2)--(1.4,-1.2);
        \draw[black,thick](.4,.2)--(1.6,-1);
        \draw[black,thick](.8,.2)--(1.8,-.8);
        \draw[black,thick](1.2,.2)--(2,-.6);
        \draw[black,thick](1.6,.2)--(2.2,-.4);
        \draw[black,thick](2,.2)--(2.4,-.2);
        \draw[black,thick](2.4,.2)--(2.6,0);
        \draw[black,thick](0,.2)--(-.2,0);
        \draw[black,thick](.4,.2)--(0,-.2);
        \draw[black,thick](.8,.2)--(.2,-.4);
        \draw[black,thick](1.2,.2)--(.4,-.6);
        \draw[black,thick](1.6,.2)--(.6,-.8);
        \draw[black,thick](2,.2)--(.8,-1);
        \draw[black,thick](2.4,.2)--(1,-1.2);
        \draw[black,thick](2.6,0)--(1.2,-1.4);
        \draw[black,thick](-.2,0)--(1.2,-1.4);
	\draw (.2,.4) node {$\cdots$};
	\draw (.8,.4) node {$\bullet$};
	\draw (1.2,.4) node {$\circ$};
	\draw (1.6,.4) node {$\circ$};
	\draw (2,.4) node {$\circ$};
	\draw (2.4,.4) node {$\circ$};
	 \draw (.8,0) node {\scriptsize{$\alpha$}};
	 \draw (1,-.2) node {\scriptsize{$\alpha$}};
	 \draw (1.2,0) node {\scriptsize{$\beta$}};
	 \draw (1.6,0) node {\scriptsize{$\beta$}};
	 \draw (2,0) node {\scriptsize{$\beta$}};
	 \draw (2.4,0) node {\scriptsize{$\beta$}};
\end{tikzpicture}
}
\hspace{1cm} \resizebox{2.8in}{!}{
 \begin{tikzpicture}
	\draw (.4,-.4) node [draw,scale=.7,diamond,fill=gray!30]{};
	\draw (.6,-.2) node [draw,scale=.7,diamond,fill=gray!30]{};
	\draw (1.2,-.4) node [draw,scale=.7,diamond,fill=gray!30]{};
	\draw (1.4,-.6) node [draw,scale=.7,diamond,fill=gray!30]{};
	\draw (1.6,-.8) node [draw,scale=.7,diamond,fill=gray!30]{};
        \draw[black,thick](0,.2)--(1.4,-1.2);
        \draw[black,thick](.4,.2)--(1.6,-1);
        \draw[black,thick](.8,.2)--(1.8,-.8);
        \draw[black,thick](1.2,.2)--(2,-.6);
        \draw[black,thick](1.6,.2)--(2.2,-.4);
        \draw[black,thick](2,.2)--(2.4,-.2);
        \draw[black,thick](2.4,.2)--(2.6,0);
        \draw[black,thick](0,.2)--(-.2,0);
        \draw[black,thick](.4,.2)--(0,-.2);
        \draw[black,thick](.8,.2)--(.2,-.4);
        \draw[black,thick](1.2,.2)--(.4,-.6);
        \draw[black,thick](1.6,.2)--(.6,-.8);
        \draw[black,thick](2,.2)--(.8,-1);
        \draw[black,thick](2.4,.2)--(1,-1.2);
        \draw[black,thick](2.6,0)--(1.2,-1.4);
        \draw[black,thick](-.2,0)--(1.2,-1.4);
	\draw (.2,.4) node {$\cdots$};
	\draw(1,.6) node {{$\curvearrowright$}};
	\draw (.8,.4) node {$\bullet$};
	\draw (1.2,.4) node {$\circ$};
	\draw (1.6,.4) node {$\bullet$};
	\draw (2,.4) node {$\circ$};
	 \draw (.8,0) node {\scriptsize{$\alpha$}};
	 \draw (1,-.2) node {\scriptsize{$\beta$}};
	 \draw (1.2,0) node {\scriptsize{$\beta$}};
	 \draw (1.6,0) node {\scriptsize{$\alpha$}};
	 \draw (2,0) node {\scriptsize{$\beta$}};
	 \draw (2.8,-.7) node {$1$};
	 \draw (2.8,-.9) node {$\mapsto$};
\end{tikzpicture}
		\begin{tikzpicture}
	\draw (.8,-.8) node [draw,scale=.7,diamond,fill=gray!30]{};
	\draw (1,-.6) node [draw,scale=.7,diamond,fill=gray!30]{};
	\draw (1.2,-.4) node [draw,scale=.7,diamond,fill=gray!30]{};
	\draw (2,-.4) node [draw,scale=.7,diamond,fill=gray!30]{};
	\draw (1.4,-.2) node [draw,scale=.7,diamond,fill=gray!30]{};
        \draw[black,thick](0,.2)--(1.4,-1.2);
        \draw[black,thick](.4,.2)--(1.6,-1);
        \draw[black,thick](.8,.2)--(1.8,-.8);
        \draw[black,thick](1.2,.2)--(2,-.6);
        \draw[black,thick](1.6,.2)--(2.2,-.4);
        \draw[black,thick](2,.2)--(2.4,-.2);
        \draw[black,thick](2.4,.2)--(2.6,0);
        \draw[black,thick](0,.2)--(-.2,0);
        \draw[black,thick](.4,.2)--(0,-.2);
        \draw[black,thick](.8,.2)--(.2,-.4);
        \draw[black,thick](1.2,.2)--(.4,-.6);
        \draw[black,thick](1.6,.2)--(.6,-.8);
        \draw[black,thick](2,.2)--(.8,-1);
        \draw[black,thick](2.4,.2)--(1,-1.2);
        \draw[black,thick](2.6,0)--(1.2,-1.4);
        \draw[black,thick](-.2,0)--(1.2,-1.4);
	\draw (.2,.4) node {$\cdots$};
	\draw (.8,.4) node {$\circ$};
	\draw (1.2,.4) node {$\bullet$};
	\draw (1.6,.4) node {$\bullet$};
	\draw (2,.4) node {$\circ$};
	 \draw (.8,0) node {\scriptsize{$\beta$}};
	 \draw (1.2,0) node {\scriptsize{$\alpha$}};
	 \draw (1.6,0) node {\scriptsize{$\alpha$}};
	 \draw (1.8,-.2) node {\scriptsize{$\beta$}};
	 \draw (2,0) node {\scriptsize{$\beta$}};
\end{tikzpicture}
}

\vspace{.3cm}
\resizebox{2.8in}{!}{
 		\begin{tikzpicture}
	\draw (.6,-.6) node [draw,scale=.7,diamond,fill=gray!30]{};
	\draw (.8,-.4) node [draw,scale=.7,diamond,fill=gray!30]{};
	\draw (1,-.2) node [draw,scale=.7,diamond,fill=gray!30]{};
	\draw (1.2,-.4) node [draw,scale=.7,diamond,fill=gray!30]{};
	\draw (1.4,-.6) node [draw,scale=.7,diamond,fill=gray!30]{};
	\draw (1.6,-.8) node [draw,scale=.7,diamond,fill=gray!30]{};
        \draw[black,thick](0,.2)--(1.4,-1.2);
        \draw[black,thick](.4,.2)--(1.6,-1);
        \draw[black,thick](.8,.2)--(1.8,-.8);
        \draw[black,thick](1.2,.2)--(2,-.6);
        \draw[black,thick](1.6,.2)--(2.2,-.4);
        \draw[black,thick](2,.2)--(2.4,-.2);
        \draw[black,thick](2.4,.2)--(2.6,0);
        \draw[black,thick](0,.2)--(-.2,0);
        \draw[black,thick](.4,.2)--(0,-.2);
        \draw[black,thick](.8,.2)--(.2,-.4);
        \draw[black,thick](1.2,.2)--(.4,-.6);
        \draw[black,thick](1.6,.2)--(.6,-.8);
        \draw[black,thick](2,.2)--(.8,-1);
        \draw[black,thick](2.4,.2)--(1,-1.2);
        \draw[black,thick](2.6,0)--(1.2,-1.4);
        \draw[black,thick](-.2,0)--(1.2,-1.4);
	\draw (.2,.4) node {$\cdots$};
	\draw(1,.6) node {{$\curvearrowleft$}};
	\draw (.8,.4) node {$\circ$};
	\draw (1.2,.4) node {$\bullet$};
	\draw (2,.4) node {$\cdots$};
	 \draw (.8,0) node {\scriptsize{$\beta$}};
	 \draw (1.2,0) node {\scriptsize{$\alpha$}};
	 \draw (2.8,-.7) node {$q$};
	\draw(-.2,.6) node {};
	 \draw (2.8,-.9) node {$\mapsto$};
\end{tikzpicture}
 \begin{tikzpicture}
	\draw (.4,-.4) node [draw,scale=.7,diamond,fill=gray!30]{};
	\draw (.6,-.2) node [draw,scale=.7,diamond,fill=gray!30]{};
	\draw (1.4,-.2) node [draw,scale=.7,diamond,fill=gray!30]{};
	\draw (1.6,-.4) node [draw,scale=.7,diamond,fill=gray!30]{};
	\draw (1.8,-.6) node [draw,scale=.7,diamond,fill=gray!30]{};
        \draw[black,thick](0,.2)--(1.4,-1.2);
        \draw[black,thick](.4,.2)--(1.6,-1);
        \draw[black,thick](.8,.2)--(1.8,-.8);
        \draw[black,thick](1.2,.2)--(2,-.6);
        \draw[black,thick](1.6,.2)--(2.2,-.4);
        \draw[black,thick](2,.2)--(2.4,-.2);
        \draw[black,thick](2.4,.2)--(2.6,0);
        \draw[black,thick](0,.2)--(-.2,0);
        \draw[black,thick](.4,.2)--(0,-.2);
        \draw[black,thick](.8,.2)--(.2,-.4);
        \draw[black,thick](1.2,.2)--(.4,-.6);
        \draw[black,thick](1.6,.2)--(.6,-.8);
        \draw[black,thick](2,.2)--(.8,-1);
        \draw[black,thick](2.4,.2)--(1,-1.2);
        \draw[black,thick](2.6,0)--(1.2,-1.4);
        \draw[black,thick](-.2,0)--(1.2,-1.4);
	\draw (.2,.4) node {$\cdots$};
	\draw (.8,.4) node {$\bullet$};
	\draw (1.2,.4) node {$\circ$};
	\draw (2,.4) node {$\cdots$};
	 \draw (.8,0) node {\scriptsize{$\alpha$}};
	 \draw (1,-.2) node {\scriptsize{$q$}};
	 \draw (1.2,0) node {\scriptsize{$\beta$}};
\end{tikzpicture}
}
%
\hspace{1cm}
\resizebox{2.8in}{!}{
 \begin{tikzpicture}
	\draw (.4,-.4) node [draw,scale=.7,diamond,fill=gray!30]{};
	\draw (.6,-.2) node [draw,scale=.7,diamond,fill=gray!30]{};
	\draw (1.4,-.2) node [draw,scale=.7,diamond,fill=gray!30]{};
	\draw (1.6,-.4) node [draw,scale=.7,diamond,fill=gray!30]{};
	\draw (1.8,-.6) node [draw,scale=.7,diamond,fill=gray!30]{};
        \draw[black,thick](0,.2)--(1.4,-1.2);
        \draw[black,thick](.4,.2)--(1.6,-1);
        \draw[black,thick](.8,.2)--(1.8,-.8);
        \draw[black,thick](1.2,.2)--(2,-.6);
        \draw[black,thick](1.6,.2)--(2.2,-.4);
        \draw[black,thick](2,.2)--(2.4,-.2);
        \draw[black,thick](2.4,.2)--(2.6,0);
        \draw[black,thick](0,.2)--(-.2,0);
        \draw[black,thick](.4,.2)--(0,-.2);
        \draw[black,thick](.8,.2)--(.2,-.4);
        \draw[black,thick](1.2,.2)--(.4,-.6);
        \draw[black,thick](1.6,.2)--(.6,-.8);
        \draw[black,thick](2,.2)--(.8,-1);
        \draw[black,thick](2.4,.2)--(1,-1.2);
        \draw[black,thick](2.6,0)--(1.2,-1.4);
        \draw[black,thick](-.2,0)--(1.2,-1.4);
	\draw (.2,.4) node {$\cdots$};
	\draw(1,.6) node {{$\curvearrowright$}};
	\draw (.8,.4) node {$\bullet$};
	\draw (1.2,.4) node {$\circ$};
	\draw (2,.4) node {$\cdots$};
	 \draw (.8,0) node {\scriptsize{$\alpha$}};
	 \draw (1,-.2) node {\scriptsize{$q$}};
	 \draw (1.2,0) node {\scriptsize{$\beta$}};
	 \draw (2.8,-.7) node {$1$};
	 \draw (2.8,-.9) node {$\mapsto$};
\end{tikzpicture}
		\begin{tikzpicture}
	\draw (.6,-.6) node [draw,scale=.7,diamond,fill=gray!30]{};
	\draw (.8,-.4) node [draw,scale=.7,diamond,fill=gray!30]{};
	\draw (1,-.2) node [draw,scale=.7,diamond,fill=gray!30]{};
	\draw (1.2,-.4) node [draw,scale=.7,diamond,fill=gray!30]{};
	\draw (1.4,-.6) node [draw,scale=.7,diamond,fill=gray!30]{};
	\draw (1.6,-.8) node [draw,scale=.7,diamond,fill=gray!30]{};
        \draw[black,thick](0,.2)--(1.4,-1.2);
        \draw[black,thick](.4,.2)--(1.6,-1);
        \draw[black,thick](.8,.2)--(1.8,-.8);
        \draw[black,thick](1.2,.2)--(2,-.6);
        \draw[black,thick](1.6,.2)--(2.2,-.4);
        \draw[black,thick](2,.2)--(2.4,-.2);
        \draw[black,thick](2.4,.2)--(2.6,0);
        \draw[black,thick](0,.2)--(-.2,0);
        \draw[black,thick](.4,.2)--(0,-.2);
        \draw[black,thick](.8,.2)--(.2,-.4);
        \draw[black,thick](1.2,.2)--(.4,-.6);
        \draw[black,thick](1.6,.2)--(.6,-.8);
        \draw[black,thick](2,.2)--(.8,-1);
        \draw[black,thick](2.4,.2)--(1,-1.2);
        \draw[black,thick](2.6,0)--(1.2,-1.4);
        \draw[black,thick](-.2,0)--(1.2,-1.4);
	\draw (.2,.4) node {$\cdots$};
	\draw (.8,.4) node {$\circ$};
	\draw (1.2,.4) node {$\bullet$};
	\draw (2,.4) node {$\cdots$};
	 \draw (.8,0) node {\scriptsize{$\beta$}};
	 \draw (1.2,0) node {\scriptsize{$\alpha$}};
\end{tikzpicture}
}
%
	\caption{Transitions in the Markov chain on tableaux.}\label{fig:5}
\end{figure}
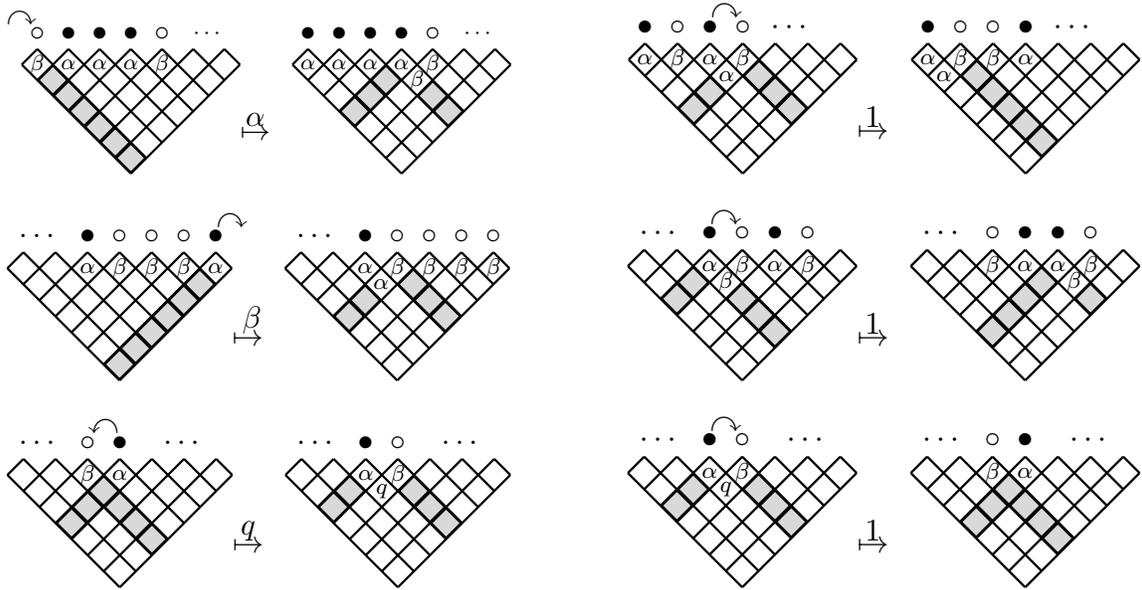

If we identify $\alpha$'s and $\beta$'s inside the tableaux
with particles and holes, then 
the chain on tableaux reveals a circulation of particles and holes
in the second row of the tableaux; this is similar to a phenomenon
observed in \cite{DS}.


\subsection{The Matrix Ansatz}
One can also  prove \cref{thm:1} using
the following \emph{Matrix Ansatz}, first introduced by 
 Derrida, Evans, Hakim and Pasquier \cite{Derrida1}.\footnote{\cite{Derrida1}
stated this result with $c=1$, but we will use
a representation with $c=\alpha \beta$ in order to 
prove \cref{thm:1}.}

	\begin{theorem}[Derrida-Evans-Hakim-Pasquier]  \label{ansatz}
Consider the ASEP with open boundaries
on a lattice of $n$ sites.
Suppose that $D$ and $E$ are matrices,  $|V\rangle$ is a column vector,
$\langle W|$ is a row vector, and $c$ is a constant,
such that:
\begin{align}
	& DE - qED = c(D+E) \label{eq:1} \\
	& \beta D|V\rangle =  c |V\rangle \label{eq:2} \\
	& \alpha \langle W|E =  c \langle W| \label{eq:3}
\end{align}
If we identify $\tau = (\tau_1,\dots,\tau_n)\in \{0,1\}^n$ with a state
(by mapping $1$ and $0$ to $\bullet$ and $\circ$, respectively), then 
	the steady state probability $\pi(\tau)$ is equal to 
\begin{equation*}
	{\pi}(\ttt) = \frac{\langle W|(\prod_{i=1}^n (\ttt_i D + (1-\ttt_i)E))|V\rangle}{\langle W|(D+E)^n|V\rangle}.
\end{equation*}
\end{theorem}

For example,  the steady state probability of 
state $\circ \bullet \circ \bullet \bullet \circ \bullet$ is 
$\frac{\langle W|EDEDDED|V\rangle}{\langle W|(D+E)^7|V\rangle}$.

We note that \cref{ansatz}  does not 
imply that a solution $D, E, |V\rangle, \langle W|$ exists nor that it is unique.  
Indeed there are multiple solutions, which in general involve 
infinite-dimensional matrices.

To prove \cref{thm:1} using the Matrix Ansatz, we let 
 $D_1=(d_{ij})$ be the 
(infinite) upper-triangular matrix with rows and columns 
indexed by $\mathbb{Z}^+$,
defined by 
$	d_{i,i+1} =\alpha$ and  $d_{ij}=0 \text{ for }j\neq i+1.$
Let 
 $E_1=(e_{ij})$ be the (infinite) lower-triangular matrix defined by 
	 $e_{ij}=0$ for $j>i$ and 
\begin{equation*}
e_{ij}=\beta^{i-j+1}(q^{j-1} {i-1 \choose j-1} + \alpha \sum_{r=0}^{j-2} {i-j+r \choose r} q^r)
	\text{ for }j\leq i.
\end{equation*}  
That is, 
\begin{equation*}
D_1=\left( \begin{array}{ccccc}
0 & \alpha & 0 & 0 & \dots \\
0 & 0 & \alpha & 0 & \dots \\
0 & 0 & 0 & \alpha & \dots \\
0 & 0 & 0 & 0 & \dots \\
\vdots & \vdots & \vdots & \vdots &
\end{array} \right)
	\text{ and }
E_1 =  \left( \begin{array}{ccccc}
\beta & 0 & 0 &   \dots \\
	\beta^2 & \beta(\alpha+q) & 0 &  \dots \\
	\beta^3 & \beta^2(\alpha+2 q) &  \beta (\alpha+\alpha q+ q^2) &
\dots \\
 \beta^4  & \beta^3 (\alpha +3 q) &
	\beta^2(\alpha+2\alpha q+3 q^2)  &  \dots \\
\vdots & \vdots & \vdots & 
\end{array} \right). 
\end{equation*}

\normalsize{
	We also define the (infinite) row and column vectors
$\langle W_1|=(1,0,0,\dots )$ and  $|V_1 \rangle= (1,1,1,\dots )^t$.
Then one can check that 
$D_1, E_1, \langle W_1|, |V_1\rangle$ satisfy  
\eqref{eq:1}, \eqref{eq:2}, and \eqref{eq:3},
with $c = \alpha \beta$.
One can also show
$D_1$ and $E_1$ are \emph{transfer matrices} whose products enumerate
$\alpha \beta$-staircase tableaux. 
For example, 
$\langle W_1|E_1D_1E_1D_1D_1E_1D_1|V_1\rangle$  enumerates the staircase tableaux
of type 
 $\circ \bullet \circ \bullet \bullet \circ \bullet$.
Now \cref{ansatz} implies \cref{thm:1}.
}

\subsection{Generalization to the five-parameter open boundary ASEP}

More generally, we would like to understand a generalized ASEP 
in which particles can both enter and exit the lattice at the left
(at rates $\alpha$, $\gamma$), and exit and enter the lattice
at the right (at rates $\beta$, $\delta$).  There is a version of the 
Matrix Ansatz for this setting \cite{Derrida1}, as well as suitable tableaux filled with 
$\alpha, \beta, \gamma$ and $\delta$'s 
(which we will simply call staircase tableaux) 
\cite{CW4}.

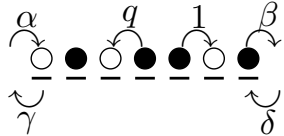
\begin{figure}[h]
\resizebox{1.7in}{!}{
\begin{tikzpicture}
        \draw[black,thick](-.1,-.2)--(.1,-.2);
        \draw[black,thick](.25,-.2)--(.45,-.2);
        \draw[black,thick](.6,-.2)--(.8,-.2);
        \draw[black,thick](.95,-.2)--(1.15,-.2);
        \draw[black,thick](1.3,-.2)--(1.5,-.2);
        \draw[black,thick](1.65,-.2)--(1.85,-.2);
        \draw[black,thick](2,-.2)--(2.2,-.2);
        \filldraw[color=black,fill=white] (0,0) circle (3pt);
        \filldraw[color=black,fill=black] (.35,0) circle (3pt);
        \filldraw[color=black,fill=white] (.7,0) circle (3pt);
        \filldraw[color=black,fill=black] (1.05,0) circle (3pt);
        \filldraw[color=black,fill=black] (1.4,0) circle (3pt);
        \filldraw[color=black,fill=white] (1.75,0) circle (3pt);
        \filldraw[color=black,fill=black] (2.1,0) circle (3pt);
	 \draw (-.15,.4) node {\small{$\alpha$}};
	 \draw (-.15,-.65) node {\small{$\gamma$}};
	  \draw(-.15,-.4) node {$\swingleft$};
	 \draw (-.15,.2) node {$\curvearrowright$};
	  \draw(.85,.2) node {$\curvearrowleft$};
	  \draw(.9,.42) node {\small{$q$}};
	  \draw(1.6,.2) node {$\curvearrowright$};
	  \draw(1.6,.42) node {\small{$1$}};
	  \draw(2.25,.2) node {$\curvearrowright$};
	  \draw(2.25,-.4) node {$\swingleft$};
	  \draw(2.3,.42) node {\small{$\beta$}};
	  \draw(2.3,-.65) node {\small{$\delta$}};
\end{tikzpicture}
}
	\caption{The (five-parameter) open boundary ASEP.} \label{fig:6}
\end{figure}

	\begin{definition} A
		\emph{staircase tableau}
	$T$ of size $n$ is a 
Young diagram of shape $(n, n-1, \dots, 2, 1)$ 
such that each box is either empty
or contains an $\alpha$, $\beta$, $\gamma$, or $\delta$ such that:
\begin{enumerate}
	\item no box in the top row is empty
	\item each box southeast of a $\beta$ or $\delta$ and in the same diagonal
 as that $\beta$ or $\delta$ is empty.
 \item each box southwest of an $\alpha$ or $\gamma$ and in the same diagonal
 as that $\alpha$ or $\gamma$ is empty.
\end{enumerate}
\end{definition}

	\begin{figure}[h]
\resizebox{4in}{!}{
\begin{tikzpicture}
        \draw[black,thick](0,.2)--(1.4,-1.2);
        \draw[black,thick](.4,.2)--(1.6,-1);
        \draw[black,thick](.8,.2)--(1.8,-.8);
        \draw[black,thick](1.2,.2)--(2,-.6);
        \draw[black,thick](1.6,.2)--(2.2,-.4);
        \draw[black,thick](2,.2)--(2.4,-.2);
        \draw[black,thick](2.4,.2)--(2.6,0);
        \draw[black,thick](0,.2)--(-.2,0);
        \draw[black,thick](.4,.2)--(0,-.2);
        \draw[black,thick](.8,.2)--(.2,-.4);
        \draw[black,thick](1.2,.2)--(.4,-.6);
        \draw[black,thick](1.6,.2)--(.6,-.8);
        \draw[black,thick](2,.2)--(.8,-1);
        \draw[black,thick](2.4,.2)--(1,-1.2);
        \draw[black,thick](2.6,0)--(1.2,-1.4);
        \draw[black,thick](-.2,0)--(1.2,-1.4);
	\draw (0,.4) node {$\circ$};
	\draw (.4,.4) node {$\bullet$};
	\draw (.8,.4) node {$\circ$};
	\draw (1.2,.4) node {$\bullet$};
	\draw (1.6,.4) node {$\bullet$};
	\draw (2,.4) node {$\circ$};
	\draw (2.4,.4) node {$\bullet$};
	 \draw (0,0) node {\scriptsize{$\gamma$}};
	 \draw (.4,0) node {\scriptsize{$\alpha$}};
	 \draw (.8,0) node {\scriptsize{$\beta$}};
	 \draw (.8,-.4) node {\scriptsize{$\gamma$}};
	 \draw (1.2,0) node {\scriptsize{$\delta$}};
	 \draw (1.2,-.8) node {\scriptsize{$\delta$}};
	 \draw (1.6,0) node {\scriptsize{$\alpha$}};
	 \draw (2,0) node {\scriptsize{$\beta$}};
	 \draw (2,-.4) node {\scriptsize{$\alpha$}};
	 \draw (2.4,0) node {\scriptsize{$\delta$}};
\end{tikzpicture}\hspace{1in}
\begin{tikzpicture}
        \draw[black,thick](0,.2)--(1.4,-1.2);
        \draw[black,thick](.4,.2)--(1.6,-1);
        \draw[black,thick](.8,.2)--(1.8,-.8);
        \draw[black,thick](1.2,.2)--(2,-.6);
        \draw[black,thick](1.6,.2)--(2.2,-.4);
        \draw[black,thick](2,.2)--(2.4,-.2);
        \draw[black,thick](2.4,.2)--(2.6,0);
        \draw[black,thick](0,.2)--(-.2,0);
        \draw[black,thick](.4,.2)--(0,-.2);
        \draw[black,thick](.8,.2)--(.2,-.4);
        \draw[black,thick](1.2,.2)--(.4,-.6);
        \draw[black,thick](1.6,.2)--(.6,-.8);
        \draw[black,thick](2,.2)--(.8,-1);
        \draw[black,thick](2.4,.2)--(1,-1.2);
        \draw[black,thick](2.6,0)--(1.2,-1.4);
        \draw[black,thick](-.2,0)--(1.2,-1.4);
	 \draw (0,0) node {\scriptsize{$\gamma$}};
	 \draw (.4,0) node {\scriptsize{$\alpha$}};
	 \draw (.4,-.4) node {\scriptsize{$q$}};
	 \draw (.6,-.2) node {\scriptsize{$q$}};
	 \draw (.6,-.6) node {\scriptsize{$q$}};
	 \draw (.8,0) node {\scriptsize{$\beta$}};
	 \draw (.8,-.4) node {\scriptsize{$\gamma$}};
	 \draw (1.2,0) node {\scriptsize{$\delta$}};
	 \draw (1.4,-.2) node {\scriptsize{$q$}};
	 \draw (1.6,-.4) node {\scriptsize{$q$}};
	 \draw (1.8,-.6) node {\scriptsize{$q$}};
	 \draw (1.2,-.8) node {\scriptsize{$\delta$}};
	 \draw (1.4,-1) node {\scriptsize{$q$}};
	 \draw (1.6,0) node {\scriptsize{$\alpha$}};
	 \draw (1.8,-.2) node {\scriptsize{$q$}};
	 \draw (2,0) node {\scriptsize{$\beta$}};
	 \draw (2,-.4) node {\scriptsize{$\alpha$}};
	 \draw (2.4,0) node {\scriptsize{$\delta$}};
\end{tikzpicture}
}
\caption{At left: a staircase tableau $T$
		of type $(\circ \bullet \circ \bullet \bullet
		\circ \bullet)$. At right: $T$
		with a $q$ in the appropriate boxes.
		We have $\wt(T)=\alpha^3 \beta^2
		\gamma^2 \delta^3 q^8$.} \label{fig:7}
\end{figure}
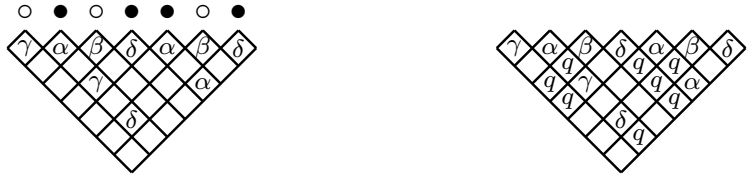

See Figure 7 for an example.
It is an exercise to verify that there are exactly $4^n n!$ staircase tableaux
of size $n$.

\begin{definition}
We call an empty box of a staircase tableau $T$
\emph{distinguished}  if either:\\
$\bullet$ its nearest neighbor on the diagonal to the northwest
is a $\delta$, or  \\
$\bullet$ its nearest neighbor on the diagonal to the northwest
is an $\alpha$ or $\gamma$, and its nearest neighbor
on the diagonal to the northeast is a $\beta$ or $\gamma$.

	After placing a $q$ in each distinguished box,
 we define the 
\emph{weight} $\wt(T)$ of $T$ to be the product of all
letters in the boxes of $T$.

The \emph{type} of $T$ is the word
obtained by reading the letters in the top row of $T$
and replacing each $\alpha$ or $\delta$
by $\bullet$, and each $\beta$ or $\gamma$ by 
$\circ$, 
see Figure 7.
\end{definition}

The following result from \cite{CWPNAS, CW4}
subsumes \cref{thm:1}.  It can be proved using
a suitable generalization
of the Matrix Ansatz.  

\begin{theorem}\label{thm:2}
Consider the ASEP with open boundaries on a lattice of $n$ sites
as in Figure 6.
Let $\tau\in \{\bullet,\circ\}^n$ be
a state.  Then the 
unnormalized steady state probability $\Psi(\tau)$ is equal to 
 $\sum_T \wt(T)$, where the sum is over the 
staircase tableaux of type $\tau$.
\end{theorem}

Remarkably, there is another solution to the Matrix Ansatz,
found earlier by Uchiyama, Sasamoto, and Wadati \cite{USW}, which makes
use of orthogonal polynomials.  More specifically, one can find a 
solution where $D$ and $E$ are tridiagonal matrices, such that the  rows of
$D+E$ encode
the three-term recurrence relation characterizing the 
\emph{Askey-Wilson polynomials};
these 
are a family of orthogonal polynomials
$p_n(x; a, b, c, d|q)$ which are at the top of the hierarchy of 
classical one-variable orthogonal polynomials (including
the others as special or limiting cases)
\cite{AW}.  

The connection of Askey-Wilson polynomials with the ASEP via \cite{USW} 
leads to applications on both sides.  On the one hand, 
it facilitates the 
computation of physical quantities in the ASEP
such as the  {phase diagram} \cite{USW};
it also leads to a relation between the ASEP
and the \emph{Askey-Wilson stochastic process} \cite{Bryc}. 
On the other hand, this connection has applications to the 
combinatorics of Askey-Wilson 
moments.  Since the 1980's there has been a great deal of 
work on the combinatorics
 of classical orthogonal polynomials (e.g. Hermite, Charlier, Laguerre)
\cite{Viennot-book,  ISV,  CKS};
the connection of staircase tableaux to ASEP, and of ASEP to Askey-Wilson polynomials, led
to the first combinatorial formula for moments of Askey-Wilson polynomials \cite{CW4, CSSW}.

Even more generally, one can study a version of the ASEP
with open boundaries in which there are different species
of particles.  This version is closely connected 
\cite{Cantini, CW5, CantiniGarbali}
to 
\emph{Koornwinder polynomials}  \cite{Koornwinder},
a family of multivariate orthogonal polynomials
which generalize Askey-Wilson polynomials.

\section{The (multispecies) ASEP on a ring}

It is also natural to consider the ASEP on
a lattice of 
 sites arranged in a ring, of which some sites
are occupied by a particle.  Each particle in the system can
jump to the next site either clockwise or counterclockwise,
provided that this site is empty.
In this model, 
the resulting stationary distribution
is always the uniform distribution.   This motivates 
considering a \emph{multispecies} generalization of the ASEP,
in which particles come with different weights, which in 
turn influence the hopping rates.



\subsection{The multispecies ASEP, multiline queues,  and Macdonald polynomials}  
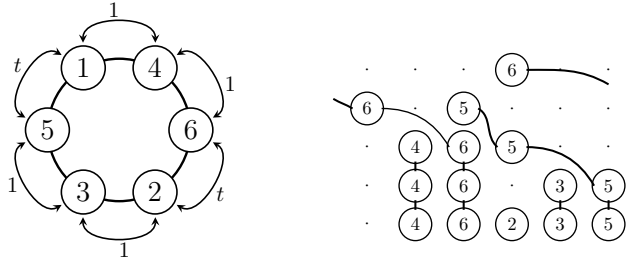
\begin{figure}[h]
\resizebox{1.5in}{!}{
	\begin{tikzpicture}[->,main node/.style={circle,draw=black, thick, fill=white,font=\Large\bfseries}]
		\draw [very thick] (0,0) circle (1.25cm);
\node[draw=none, minimum size=2.5cm, regular polygon, regular polygon sides=6] (s) {};
\node[draw=none, minimum size=3.3cm, regular polygon, regular polygon sides=24] (b) {};
\node[draw=none, minimum size=4.2cm, regular polygon, regular polygon sides=72] (outer) {};
                \node[shift=(outer.corner 2)] {$1$};
                \node[shift=(outer.corner 13)] {$t$};
                \node[shift=(outer.corner 25)] {$1$};
                \node[shift=(outer.corner 38)] {$1$};
                \node[shift=(outer.corner 49)] {$t$};
                \node[shift=(outer.corner 60)] {$1$};
		\node[main node][shift=(s.corner 1)]
		{$4$};
		\node[main node][shift=(s.corner 2)] {$1$};
		\node[main node][shift=(s.corner 3)] {$5$};
		\node[main node][shift=(s.corner 4)] {$3$};
		\node[main node][shift=(s.corner 5)] {$2$};
		\node[main node][shift=(s.corner 6)] {$6$};
\draw[thick, >=stealth, <->](b.corner 24) to [out=90, in=90] (b.corner 3);
\draw[thick, >=stealth, <->](b.corner 4) to [out=150, in=150] (b.corner 7);
\draw[thick, >=stealth, <->](b.corner 8) to [out=210, in=210] (b.corner 11);
\draw[thick, >=stealth, <->](b.corner 12) to [out=270, in=270] (b.corner 15);
\draw[thick, >=stealth, <->](b.corner 16) to [out=330, in=330] (b.corner 19);
\draw[thick, >=stealth, <->](b.corner 20) to [out=30, in=30] (b.corner 23);
	\end{tikzpicture}\hspace{.5cm}
		}\hspace{.5cm}
\resizebox{1.6in}{!}{
\begin{tikzpicture}[main node/.style={circle,draw=black, very thick, fill=white,font=\Large\bfseries}]
     \draw (0,0) node {$\cdot$};
     \draw (0,.8) node {$\cdot$};
     \draw (0,1.6) node {$\cdot$};
    \node[circle,draw, thick] at (0,2.4){$6$};
     \draw (0,3.2) node {$\cdot$};
	\node[circle,draw, thick] at (1,0){$4$};
	\node[circle,draw, thick] at (1,.8){$4$};
	\node[circle,draw, thick] at (1,1.6){$4$};
     \draw (1,2.4) node {$\cdot$};
     \draw (1,3.2) node {$\cdot$};
	\node[circle,draw, thick] at (2,0){$6$};
	\node[circle,draw, thick] at (2,.8){$6$};
	\node[circle,draw, thick] at (2,1.6){$6$};
	\node[circle,draw, thick] at (2,2.4){$5$};
     \draw (2,3.2) node {$\cdot$};
	\node[circle,draw, thick] at (3,0){$2$};
     \draw (3,.8) node {$\cdot$};
	\node[circle,draw, thick] at (3,1.6){$5$};
     \draw (3,2.4) node {$\cdot$};
	\node[circle,draw, thick] at (3,3.2){$6$};
	\node[circle,draw, thick] at (4,0){$3$};
	\node[circle,draw, thick] at (4,.8){$3$};
     \draw (4,1.6) node {$\cdot$};
     \draw (4,2.4) node {$\cdot$};
     \draw (4,3.2) node {$\cdot$};
	\node[circle,draw, thick] at (5,0){$5$};
	\node[circle,draw, thick] at (5,.8){$5$};
     \draw (5,1.6) node {$\cdot$};
     \draw (5,2.4) node {$\cdot$};
     \draw (5,3.2) node {$\cdot$};
        \draw[black,very thick](1,1.1)--(1,1.3);
     \draw (0,-.7) node {};
        \draw[black,very thick](1,.3)--(1,.5);
        \draw[black,very thick](2,1.1)--(2,1.3);
        \draw[black,very thick](2,.3)--(2,.5);
        \draw[black,very thick](5,.3)--(5,.5);
        \draw[black,very thick](4,.3)--(4,.5);
	\draw[thick, >=stealth](.3,2.4) to [out=0, in=120] (1.7,1.6);
	\draw[very thick, >=stealth](2.3,2.4) to [out=-30, in=140] (2.7,1.6);
	\draw[very thick, >=stealth](3.3,1.6) to [out=0, in=120] (4.7,.8);
	\draw[very thick, >=stealth](3.3,3.2) to [out=0, in=150] (5,2.9);
	\draw[very thick, >=stealth](-.7,2.6) to [out=-30, in=150] (-.3,2.4);
	\end{tikzpicture}
}
		\caption{The mASEP on a ring
		with $\lambda = (6,5,4,3,2,1)$, and a multiline queue of 
		type $(1,4,6,2,3,5)$.}\label{fig:8}
\end{figure}

In the multispecies ASEP (mASEP) on a ring, two neighboring
particles exchange places at rates $1$ or $t$, depending
on whether the heavier particle is clockwise or 
counterclockwise from the lighter one.

	\begin{definition}
 Let $t$ be a constant such that $0 \leq t \leq 1$, and 
let $\lambda = \lambda_1 \geq \lambda_2 \geq \dots 
\geq \lambda_n \geq 0$
be a partition. 
Let $B_n(\lambda)$
 be the set of all words of length $n$
obtained by permuting the parts of $\lambda$.
The \emph{multispecies ASEP}
	on a ring
        is the Markov chain on $B_n(\lambda)$
        with transition probabilities:
\begin{itemize}
	\item If $\mu = (\mu_1,\dots,\mu_n)$ and $\nu$ are in $B_n(\lambda)$,
		and $\nu$ is obtained from $\mu$ by swapping 
		$\mu_i$ and $\mu_{i+1}$ for some $i$ (indices considered modulo $n$),
then
		$\Pr(\mu \to \nu) = \frac{t}{n}$ if $\mu_i>\mu_{i+1}$ and 
		$\Pr(\mu \to \nu) = \frac{1}{n}$ if $\mu_i<\mu_{i+1}$.
\item Otherwise $\Pr(\mu \to \nu) = 0$ for $\nu \neq \mu$ and
	$\Pr(\mu \to \mu) = 1-\sum_{\nu \neq \mu} \Pr(\mu \to \nu)$.
\end{itemize}
We think of the parts of $\lambda$ 
as representing various types of
particles of different weights.
	\end{definition}
	
As before, one would like to find an expression for each steady
state probability  as a manifestly positive
sum over some set of combinatorial objects.  
One may give such a
formula in terms of Ferrari-Martin's
\emph{multiline queues} shown in Figure 8,
see \cite{Martin, CMW2}.

One fascinating aspect %
of the multispecies
ASEP on a ring is its close relation \cite{CGW} 
to \emph{Macdonald polynomials} $P_{\lambda}(x_1,\dots,x_n; q,t)$ \cite{Macdonald},
a remarkable family of polynomials that generalize 
Schur polynomials, Hall-Littlewood polynomials, and Jack 
polynomials.   The next result 
follows  from \cite{CGW} and \cite{CMW2}.

\begin{theorem}\label{thm:Mac}
Let $\mu \in B_n(\lambda)$ be a state of the 
mASEP on a ring.  Then the steady state probability 
$\pi(\mu)$ 
	is $$\pi(\mu)=\frac{\Psi(\mu)}{Z_{\lambda}},$$
	where $\Psi({\mu})$ 
is obtained from a \emph{permuted basement
	Macdonald polynomial} and 
	 $Z_{\lambda}$
	is obtained from the \emph{Macdonald polynomial}
$P_{\lambda}$ by specializing
$q=1$ and $x_1=x_2=\dots=x_n=1$.
\end{theorem}
	
The following table shows the  probabilities of the mASEP when 
$\lambda=(4,3,2,1)$.  Note that because of the circular symmetry in the mASEP,
e.g. $\pi(1,2,3,4)=\pi(2,3,4,1)=\pi(3,4,1,2)=\pi(4,1,2,3)$, it suffices to list
the probabilities for the states $w$ with $w_1=1$.

\begin{table}[h]
\begin{center}
\begin{tabular}{|c c| }
    \hline
	State $w$ & Unnormalized probability $\Psi(w)$\\
    \hline 
    1234 & $9t^3+7t^2+7t+1$\\
	1243 & $3(t^3+3t^2+3t+1)$\\
	1324 & $3t^3+11t^2+5t+5$\\
	1342 & $3(t^3+3t^2+3t+1)$\\
	1423 & $5t^3+5t^2+11t+3$\\
	1432 & $t^3+7t^2+7t+9$\\ 
    \hline
    \end{tabular}
\end{center}
	\caption{Probabilities for the 
	mASEP when $\lambda=(4,3,2,1)$.}\label{table:2}
\end{table}
\normalsize{

In light of 
\cref{thm:Mac} and the connection to multiline queues, 
it is natural to ask if one can give a formula for Macdonald
polynomials in terms of multiline queues.   This is indeed possible,
see \cite{CMW2} for details.
}

We remark that there is a family of 
Macdonald polynomials associated to any affine root system;
the ``ordinary'' Macdonald polynomials discussed in this section
are those of type $\tilde{A}$.  It is interesting that they
are related to particles hopping on a ring (which 
resembles
an affine $A$ Dynkin diagram).  Meanwhile, the Koornwinder
polynomials from the previous section
are the 
Macdonald 
 polynomials attached to the non-reduced affine root system
 of type
$\tilde{C}_n^{\vee}$. 
It is interesting that they are related to particles hopping
on a line with open boundaries (which 
 resembles a Dynkin diagram of type 
$\tilde{C}_n^{\vee}$).

We note that there
are other connections between probability and Macdonald polynomials, including 
 \emph{Macdonald processes} \cite{BorodinCorwin}, and
a Markov chain on partitions whose 
eigenfunctions are coefficients of Macdonald polynomials \cite{Diaconis}.
There is also a variation of the exclusion process called the 
\emph{multispecies zero range process}, whose stationary distribution 
is related to 
\emph{modified Macdonald polynomials} 
\cite{AMM}.

\subsection{The inhomogeneous TASEP, multiline queues,  
and Schubert polynomials}  

Another 
multispecies generalization 
of the exclusion process
on a ring is the \emph{inhomogeneous totally asymmetric exclusion 
process} (TASEP).  In this model, two adjacent particles with weights
$i$ and $j$ with $i<j$
can swap places only if the heavier one is clockwise
of the lighter one, and in this case, they exchange places at
rate $x_i-y_j$, 
 see Figure 9.
\begin{figure}[h]
\resizebox{1.6in}{!}{
	\begin{tikzpicture}[->,main node/.style={circle,draw=black, thick, fill=white,font=\Large\bfseries}]
		\draw [very thick] (0,0) circle (1.25cm);
\node[draw=none, minimum size=2.5cm, regular polygon, regular polygon sides=6] (s) {};
\node[draw=none, minimum size=3.3cm, regular polygon, regular polygon sides=24] (b) {};
\node[draw=none, minimum size=4.2cm, regular polygon, regular polygon sides=72] (outer) {};
                \node[shift=(outer.corner 2)] {$x_1-y_4$};
                \node[shift=(outer.corner 25)] {$x_3-y_5$};
                \node[shift=(outer.corner 38)] {$x_2-y_3$};
                \node[shift=(outer.corner 60)] {$x_4-y_6$};
		\node[main node][shift=(s.corner 1)]
		{$4$};
		\node[main node][shift=(s.corner 2)] {$1$};
		\node[main node][shift=(s.corner 3)] {$5$};
		\node[main node][shift=(s.corner 4)] {$3$};
		\node[main node][shift=(s.corner 5)] {$2$};
		\node[main node][shift=(s.corner 6)] {$6$};
		\draw[thick, >=stealth, <->](b.corner 24) to [out=90, in=90] (b.corner 3);
		\draw[thick, >=stealth, <->](b.corner 8) to [out=210, in=210] (b.corner 11);
		\draw[thick, >=stealth, <->](b.corner 12) to [out=270, in=270] (b.corner 15);
		\draw[thick, >=stealth, <->](b.corner 20) to [out=30, in=30] (b.corner 23);
\end{tikzpicture}
        }
		\caption{The 
		 inhomogeneous multispecies TASEP on a ring,
		with $\lambda = (6,5,4,3,2,1)$.}\label{fig:9}
\end{figure}
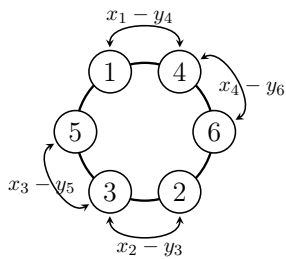

	\begin{definition}
		Let $x_1,\dots,x_n$ and $y_1,\dots,y_n$ be constants
such that $0<x_i-y_j\leq 1$ for all $i,j$, 
and let $\lambda = \lambda_1 \geq \lambda_2 \geq \dots 
\geq \lambda_n \geq 0$
be a partition. 
Let ${B}_n(\lambda)$
 be the set of all words of length $n$
obtained by permuting the parts of $\lambda$.
The \emph{inhomogeneous TASEP}
on a ring
is the Markov chain on ${B}_n(\lambda)$
      with transition probabilities:
\begin{itemize}
\item If $\mu = (\mu_1,\dots,\mu_n)$ and $\nu$ 
	are in ${B}_n(\lambda)$,
	and $\nu$ is obtained from $\mu$ by swapping
	$\mu_i$ and $\mu_{i+1}$ for some $i$
	(indices considered mod $n$), 
then
		$\Pr(\mu \to \nu) = \frac{x_{\mu_i}-y_{\mu_j}}{n}$ if $\mu_i<\mu_{i+1}$.
\item Otherwise $\Pr(\mu \to \nu) = 0$ for $\nu \neq \mu$ and
	$\Pr(\mu \to \mu) = 1-\sum_{\nu \neq \mu} \Pr(\mu \to \nu)$.
\end{itemize}
\end{definition}

When $y_i=0$ for all $i$, there is a 
formula for the stationary distribution of the 
inhomogeneous TASEP in terms of multiline queues; this can 
be proved using a version of the Matrix Ansatz
\cite{AritaMallick}.

Recall that the mASEP on a ring is closely connected to 
Macdonald polynomials.  Curiously, when $y_i=0$ the inhomogeneous TASEP on a ring
is related to \emph{Schubert polynomials}, a family of polynomials which give
polynomial representatives for the Schubert classes in the 
cohomology ring of  the 
complete flag variety.  For example, many (unnormalized) steady 
state probabilities are equal to products of Schubert polynomials
\cite{Cantini2, KW}, and all of them are conjecturally  
positive sums of Schubert polynomials \cite{LW}.

Given $w=(w_1,\dots,w_n)$ a permutation in the symmetric group $S_n$ and $p=(p_1,\dots,p_m)\in S_m$ with $m<n$, we say that $w$ \emph{contains} $p$ if $w$ has a 
subsequence of length $m$ whose letters are in the same relative order as those of $p$.
For example, the permutation $({\mathbf{3}}, 2, \mathbf{6}, 5, \mathbf{1} , \mathbf{4})$ contains the pattern
$(2,4,1,3)$ because its letters $3,6,1,4$ have the same relative 
order as those of $(2,4,1,3)$.
If $w$ does not contain $p$ we say that $w$ \emph{avoids} $p$.
We say that $w\in S_n$ is \emph{evil-avoiding} if $w$ avoids the patterns
$(2,4,1,3), (4,1,3,2), (4,2,1,3)$ and $(3,2,1,4)$.\footnote{We call these permutations \emph{evil-avoiding} because if one
replaces $i$ by $1$, $e$ by $2$, $l$ by $3$, and $v$ by $4$,
then \emph{evil} and its anagrams \emph{vile, veil} and \emph{leiv} become the four patterns $2413, 4132, 4213$ and $3214$.  (Leiv is a name of Norwegian origin meaning ``heir.'')}

We have the following result, see \cite{KW} for details.

\begin{theorem}\label{thm:KW}
Let $\lambda=(n,n-1,\dots,1)$ so that the inhomogeneous TASEP can be viewed
as a Markov chain on the $n!$ permutations of the set $\{1,2,\dots,n\}$.
Let $w\in S_n$ be a permutation with $w_1=1$ which is \emph{evil-avoiding}, and let
$k$ be the number of descents of $w^{-1}$.
	Then the steady state probability $\pi(w)$ 
	equals 
	$$\pi(w) = \frac{\Psi(w)}{Z_n},$$
	where $\Psi(w)$ is a monomial in $x_1,\dots,x_{n-1}$ 
times a product of $k$ Schubert polynomials,
and $Z_n = \prod_{i=1}^n h_{n-i}(x_1,x_2,\dots,x_{i-1},x_i,x_i)$
with $h_i$ the complete homogeneous symmetric polynomial.
\end{theorem}

The following table shows the  probabilities of the inhomogeneous 
TASEP when 
$\lambda=(4,3,2,1)$. 

\begin{table}[h]
\begin{tabular}{|c c| }
    \hline
	State $w$ & Unnormalized probability $\Psi(w)$\\
    \hline 
    1234 & $x_1^3  x_2$\\
	1243 & $x_1^2(x_1 x_2 + x_1 x_3+x_2 x_3) = x_1^2 \Sym_{1342}$\\
	1324 & $x_1(x_1^2 x_2 + x_1 x_2^2+x_1^2 x_3 + x_1 x_2 x_3 + x_2^2 x_3) =                  x_1 \Sym_{1432}$\\
	1342 & $x_1 x_2 (x_1^2+x_1 x_2 + x_2^2) = x_1 x_2 \Sym_{1423}$\\
	1423 & $x_1^2 x_2 (x_1+x_2+x_3) = x_1^2  x_2 \Sym_{1243}$\\
	1432 & $(x_1^2+x_1 x_2+x_2^2)(x_1 x_2+x_1 x_3+x_2 x_3) = \Sym_{1423} \Sym_{1342}$\\ 
    \hline
    \end{tabular}
	\caption{Probabilities for the 
	inhomogeneous TASEP when $\lambda=(4,3,2,1)$.}\label{table:3}
\end{table}

For general $y_i$, there is a version of \cref{thm:KW}
involving \emph{double} Schubert polynomials \cite{KW}.






Very often, beautiful combinatorial properties 
go hand-in-hand with \emph{integrability} of a model.
  While this topic goes beyond the scope of 
this article, the reader can learn about 
integrability and the exclusion process from \cite{Cantini2, Crampe},
or more generally
 about \emph{integrable probability} from \cite{Borodin}.

\newpage
\section{Positivity in Markov chains}

The reader may at this point wonder how general is the phenomenon of positivity in 
Markov chains?
That is, how often can one express the steady state probabilities of a Markov chain 
in terms of polynomials with all coefficients positive
(ideally as a sum over combinatorial objects)?

In some sense, the answer to this question is \emph{all the time}: 
 the \emph{Markov Chain Tree Theorem} gives a formula for the stationary 
distribution of a finite-state irreducible Markov chain as a positive
sum indexed by rooted trees
of the state diagram.  However, the number of terms of this formula grows fast
very quickly! (By Cayley's formula,
the complete graph on $n$ vertices has $n^{n-2}$ spanning trees.)
Moreover, for many Markov chains, there is a 
common factor which can be removed from the above formula for the 
stationary distribution, resulting in a more
compact formula.  Sometimes the more compact formula involves polynomials
with negative coefficients.

Let $G$ be the \emph{state diagram} of a finite-state irreducible Markov chain whose
set of states is $V$.
That is, $G$ is a weighted directed graph with vertices $V$,  with an edge $e$ from
$i$ and $j$  weighted $\Pr(e):=\Pr(i,j)$ whenever the probability $\Pr(i,j)$ of going
from state $i$ to $j$ is positive.
We call a connected subgraph $T$ a \emph{spanning tree} rooted at $r$ 
if $T$ includes every vertex of $V$, $T$ has no cycle, and all edges of $T$
point towards the root $r$.  (Irreducibility of the Markov chain implies that for 
each vertex $r$, there is a spanning tree rooted at $r$.)
Given a spanning tree $T$, we define its \emph{weight} as
$\wt(T):=\prod_{e\in T} \Pr(e).$

\begin{theorem}[Markov Chain Tree Theorem] \label{thm:Markov}
The stationary distribution of a finite-state irreducible Markov chain is proportional
	to the measure that assigns the state $\tau$ the 
	 ``unnormalized probability''
	$$\Psi(\tau):= \sum_{\rt(T)=\tau} \wt(T).$$
	That is, the steady state probability $\pi(\tau)$ equals
	$\pi(\tau) = \frac{\Psi(\tau)}{\mathbf{Z}},$
	where $\mathbf{Z} = \sum_{\tau} \Psi(\tau)$.
\end{theorem}
\cref{thm:Markov} first appeared in \cite{Hill} 
and was proved for general Markov chains in \cite{Leighton}.
It has by now many proofs, one of which involves lifting the Markov chain to a chain on the trees
themselves; the result then follows from 
\emph{Kirchhoff's Matrix Tree Theorem}.  
See
\cite{Anantharam}, 
\cite{Lyons}, 
\cite{Pitman},
and references therein.

\begin{example}\label{ex:2}
Consider the Markov chain with five states $1,\dots,5$,
whose transition matrix is as follows:
\begin{equation}
\begin{bmatrix}
	\frac{2-q}{3} & 0 & \frac{1}{3} & \frac{q}{3} & 0\\
	0 & \frac{2}{3} & 0 & 0 & \frac{1}{3}\\
	\frac{q}{3} & \frac{1}{3} & \frac{1-q}{3} & 0 & \frac{1}{3}\\
	\frac{1}{3} & \frac{q}{3} & 0 & \frac{2-2q}{3} & \frac{q}{3} \\
	0 & 0 & \frac{q}{3} & \frac{1}{3} & \frac{2-q}{3}
\end{bmatrix}
\end{equation}

	The state diagram  is shown in Figure 10.
	(We have omitted the factors of $\frac{1}{3}$
	from each 
	transition probability 
	as they do not affect
	the eigenvector of the transition matrix).  
	We also omitted the loops 
	at each state.

	\begin{figure}[h]
\resizebox{3.3in}{!}{
\begin{tikzpicture}[main node/.style={circle,draw=black, thick, fill=white,font=\bfseries}]
	\node[circle,draw, thick] at (1,1.6){$1$};
	\node[circle,draw, thick] at (1,0){$4$};
	\node[circle,draw, thick] at (2,.8){$5$};
	\node[circle,draw, thick] at (3,0){$2$};
	\node[circle,draw, thick] at (3,1.6){$3$};
	\draw[>=stealth, ->,black,thick](3, 1.3)--(3,0.4);
	 \draw (3.2,.85) node {\scriptsize{$1$}};
	\draw[>=stealth, ->,black,thick](1.3,0)--(2.6,0);
	 \draw (1.95,-.2) node {\scriptsize{$q$}};
	\draw[>=stealth, ->,black,thick](2.7,.3)--(2.3,.6);
	 \draw (2.6,.55) node {\scriptsize{$1$}};
	\draw[>=stealth, ->, black, thick](2.7,1.4) to [out=240, in=30] (2.3,1);
	 \draw (2.6,1.1) node {\scriptsize{$1$}};
	\draw[>=stealth, ->, black, thick](2.25,1.05) to [out=90, in=200] (2.65,1.45);
	 \draw (2.25,1.3) node {\scriptsize{$q$}};
	\draw[>=stealth, ->, black, thick](1.7,0.6) to [out=240, in=30] (1.3,.2);
	 \draw (1.6,.3) node {\scriptsize{$1$}};
	\draw[>=stealth, ->, black, thick](1.25,0.25) to [out=90, in=200] (1.65,0.65);
	 \draw (1.35,.6) node {\scriptsize{$q$}};
	\draw[>=stealth, ->, black, thick](.95,0.3) to [out=100, in=250] (.95,1.3);
	 \draw (.8,.8) node {\scriptsize{$1$}};
	\draw[>=stealth, ->, black, thick](1.05,1.3) to [out=-70, in=80] (1.05,.3);
	 \draw (1.2,.8) node {\scriptsize{$q$}};
	\draw[>=stealth, ->, black, thick](1.3,1.6) to [out=-20, in=190] (2.6,1.6);
	 \draw (2,1.4) node {\scriptsize{$1$}};
	\draw[>=stealth, ->, black, thick](2.6,1.7) to [out=170, in=10] (1.3,1.7);
	 \draw (2,1.9) node {\scriptsize{$q$}};
\end{tikzpicture}\hspace{1.8cm}
\begin{tikzpicture}[main node/.style={circle,draw=black, thick, fill=white,font=\Large\bfseries}]
	\node[circle,draw, thick] at (1,1.6){$1$};
	\node[circle,draw, thick] at (1,0){$4$};
	\node[circle,draw, thick] at (2,.8){$5$};
	\node[circle,draw, thick] at (3,0){$2$};
	\node[circle,draw, thick] at (3,1.6){$3$};
	\draw[>=stealth, ->,black,thick](2.7,.3)--(2.3,.6);
	 \draw (2.6,.55) node {\scriptsize{$1$}};
	\draw[>=stealth, ->, black, thick](2.25,1.05) to [out=90, in=200] (2.65,1.45);
	 \draw (2.25,1.3) node {\scriptsize{$q$}};
	\draw[>=stealth, ->, black, thick](1.25,0.25) to [out=90, in=200] (1.65,0.65);
	 \draw (1.35,.6) node {\scriptsize{$q$}};
	\draw[>=stealth, ->, black, thick](2.6,1.7) to [out=170, in=10] (1.3,1.7);
	 \draw (2,1.9) node {\scriptsize{$q$}};
	\end{tikzpicture}
}
\caption{The state diagram
 from \cref{ex:2}, plus a spanning tree rooted at $1$ with weight $q^3$.}
 \label{fig:10}
\end{figure}
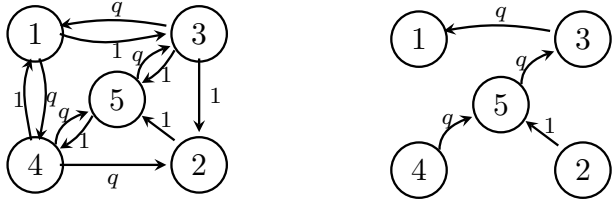

If one applies \cref{thm:Markov}, one finds e.g. that there are 
six spanning trees rooted at state $1$, with weights
$q^3, q^3, q^2, q, 1,$ and $1$.  Adding up these contributions gives 
$\Psi(1) = 2q^3+q^2+q+2$.  Computing the spanning trees
rooted at the other states gives rise to the
unnormalized probabilities $\Psi(\tau)$ for the
stationary distribution shown in \cref{table:4}.
\begin{table}[h]
\begin{center}
\begin{tabular}{|c c| }
    \hline
	State $\tau$ & Unnormalized probability $\Psi(\tau)$\\
    \hline 
	$1$ & $2q^3+q^2+q+2$\\
	$2$ & $q^4+3q^3+4q^2+3q+1$\\
	$3$  & $2q^3+2q^2+q+1$\\
	$4$  & $q^3+q^2+2q+2$\\
	$5$ & $2q^3+4q^2+4q+2$\\
    \hline
    \end{tabular}
\end{center}
	\caption{Unnormalized probabilities for the Markov
	chain from \cref{ex:2} as given by the Markov Chain Tree Theorem.}
	\label{table:4}
\end{table}

Note that the unnormalized probabilities from \cref{table:4} share a common factor of $(q+1)$.
Dividing by this common factor gives the  (more compact)
unnormalized probabilities $\overline{\Psi}(\tau)$ shown in 
\cref{table:5}.
\begin{table}[h]
\begin{center}
\begin{tabular}{|c c| }
    \hline
	State $\tau$ & Unnormalized probability $\overline{\Psi}(\tau)$\\
    \hline 
	$1$ & $2q^2-q+2$\\
	$2$ & $q^3+2q^2+2q+1$\\
	$3$  & $2q^2+1$\\
	$4$  & $q^2+2$\\
	$5$ & $2q^2+2q+2$\\
    \hline
    \end{tabular}
\end{center}
	\caption{Unnormalized probabilities from 
	\cref{table:4}
	 after dividing by the joint 
	common factor $(q+1)$.}
	\label{table:5}
\end{table}

We see that when we write the stationary distribution in 
	``lowest terms,'' we obtain a vector of polynomials which do \emph{not} have 
	only nonnegative coefficients.
\end{example}

This example motivates the following definitions.

\begin{definition}
Consider a measure 
$(\Psi_1,\dots,\Psi_n)$ 
on the 
set $\{1,2,\dots,n\}$ in which each component 
	$\Psi_i(q_1,\dots,q_N)$ is a 
	polynomial in $\Z[q_1,\dots,q_N]$.\footnote{We don't require 
	that $\sum_i \Psi_i = 1$; to obtain a probability distribution
	we can just divide each term by $Z:=\sum_i \Psi_i$.}
We say  the formula 
$(\Psi_1,\dots,\Psi_n)$ 
is \emph{manifestly positive} if 
all coefficients
of 
	$\Psi_i$ are positive for all $i$.
And we say
$(\Psi_1,\dots,\Psi_n)$ 
  is \emph{compact} 
	if there is no polynomial $\phi(q_1,\dots,q_N)\neq 1$ which divides  all the 
	$\Psi_i$.
\end{definition}

\cref{thm:Markov} shows that every finite-state Markov chain has 
a manifestly positive formula for the stationary distribution.
Meanwhile, \cref{ex:2} shows that in general this formula 
is not compact, and 
that 
there are Markov chains whose
 compact formula for the stationary distribution  is 
not manifestly positive.

In light of \cref{thm:Markov},
it is interesting to revisit 
e.g. the stationary distribution
of 
the open boundary ASEP 
with parameters
$\alpha$, $\beta$, and $q$.  
One can use \cref{thm:1} to express the components 
$\Psi_{\tab}(\tau)$ of the stationary measure
as a sum over the tableaux of type $\tau$.
On the other hand, one can use \cref{thm:Markov}
to express the components $\Psi_{\tree}(\tau)$ of the stationary measure
as a sum over spanning trees rooted at $\tau$ 
of the state diagram.
Both $\Psi_{\tab}(\tau)$ and $\Psi_{\tree}(\tau)$ are polynomials
in $\alpha, \beta, q$ with positive coefficients; however,
the former is compact, and has many fewer terms than the latter.
Because the stationary measure is unique (up to an overall scalar), 
for each $n$ there 
is a polynomial $Q_n(\alpha,\beta,q)$ such that 
$\frac{\Psi_{\tree}(\tau)}{\Psi_{\tab}(\tau)}=Q_n(\alpha, \beta,q)$
for all $\tau\in \{0,1\}^n$.  The number of terms in $Q_n$ appears in 
	\cref{table:6}.

\begin{table}[h]
\begin{center}
\begin{tabular}{|c c | }
    \hline
	 $n$ & $Q_n(1, 1, 1)$ \\
    \hline 
	$2$ & $1$\\
	$3$ & $4$\\
	$4$ & 
 $840$\\
	$5$ & 
 $2285015040 $\\
	$6$ & 
$11335132600511975880768000 $\\
    \hline
    \end{tabular}
\end{center}
	\caption{The ratio of the numbers of terms between the Markov Chain Tree theorem
	formula and the staircase tableaux formula for the stationary
	distribution of the (three-parameter) open boundary ASEP on 
	a lattice of $n$ sites.}
	\label{table:6}
\end{table}

\normalsize{

It would be interesting to reprove e.g. \cref{thm:1} 
using the Markov Chain Tree Theorem.

We note that the analysis of the ASEP and its variants would be easier
if these Markov chains were \emph{reversible}; in general 
they are not (except for special cases of the parameters).  Nevertheless
there has been progress on the mixing time of the ASEP,
see \cite{Nestoridi} and references therein.

Besides the ASEP,  
there are other interesting Markov chains
arising in statistical mechanics whose stationary distributions admit 
manifestly 
positive formulas as sums over combinatorial objects
(which are often compact).
These include 
the \emph{Razumov-Stroganov
correspondence} 
\cite{PDF, DGR, CS}, 
the \emph{Tsetlin library}
\cite{Tsetlin, Hendricks}, 
and many
other models of interacting particles \cite{AyyerMallick, AyyerNadeau}, see also
\cite{AyyerSlides}.

}

\subsection*{Acknowledgements} 
I am grateful to 
Sylvie Corteel, Olya Mandelshtam, and Donghyun Kim for some very 
stimulating and enjoyable collaborations.
I would also like to thank 
Arvind Ayyer, Alexei Borodin, Ivan Corwin, Jim Pitman, and Jim Propp for 
their comments on earlier drafts of this manuscript, which greatly 
improved the exposition.
This work was partially supported by the National Science Foundation
No.\ DMS-1854316 and No.\ DMS-1854512.



\bibliographystyle{alpha}
\bibliography{bibliography}

\def\cprime{$'$} \def\cprime{$'$}
\begin{thebibliography}{CGdGW16}

\bibitem[AM10]{AyyerMallick}
Arvind Ayyer and Kirone Mallick.
\newblock Exact results for an asymmetric annihilation process with open
  boundaries.
\newblock {\em J. Phys. A}, 43(4):045003, 22, 2010.

\bibitem[AM13]{AritaMallick}
Chikashi Arita and Kirone Mallick.
\newblock Matrix product solution of an inhomogeneous multi-species {TASEP}.
\newblock {\em J. Phys. A}, 46(8):085002, 11, 2013.

\bibitem[AMM20]{AMM}
Arvind Ayyer, Olya Mandelshtam, and James Martin.
\newblock Modified {M}acdonald polynomials and the multispecies zero range
  process: I, 2020.
\newblock Preprint, \texttt{arXiv:2011.06117}.

\bibitem[AN21]{AyyerNadeau}
Arvind Ayyer and Philippe Nadeau.
\newblock Combinatorics of a disordered two-species asep on a torus, 2021.
\newblock Preprint, \texttt{arXiv:2104.02448}.

\bibitem[AT90]{Anantharam}
Venkat Anantharam and Pantelis Tsoucas.
\newblock Stochastic concavity of throughput in series of queues with finite
  buffers.
\newblock {\em Adv. in Appl. Probab.}, 22(3):761--763, 1990.

\bibitem[AW85]{AW}
Richard Askey and James Wilson.
\newblock Some basic hypergeometric orthogonal polynomials that generalize
  {J}acobi polynomials.
\newblock {\em Mem. Amer. Math. Soc.}, 54(319):iv+55, 1985.

\bibitem[Ayy22]{AyyerSlides}
Arvind Ayyer.
\newblock Interacting particle systems and symmetric functions, 2022.
\newblock Talk at FPSAC, \texttt{https://fpsac2021.math.biu.ac.il \\
  /wp-content/uploads/2022/01/Ayyer.pdf}.

\bibitem[BC14]{BorodinCorwin}
A.~Borodin and I.~Corwin.
\newblock Macdonald processes.
\newblock In {\em X{VII}th {I}nternational {C}ongress on {M}athematical
  {P}hysics}, pages 292--316. World Sci. Publ., Hackensack, NJ, 2014.

\bibitem[BG97]{BG}
Lorenzo Bertini and Giambattista Giacomin.
\newblock Stochastic {B}urgers and {KPZ} equations from particle systems.
\newblock {\em Comm. Math. Phys.}, 183(3):571--607, 1997.

\bibitem[BG16]{Borodin}
Alexei Borodin and Vadim Gorin.
\newblock Lectures on integrable probability.
\newblock In {\em Probability and statistical physics in {S}t. {P}etersburg},
  volume~91 of {\em Proc. Sympos. Pure Math.}, pages 155--214. Amer. Math.
  Soc., Providence, RI, 2016.

\bibitem[BP18]{BorodinPetrov}
Alexei Borodin and Leonid Petrov.
\newblock Higher spin six vertex model and symmetric rational functions.
\newblock {\em Selecta Math. (N.S.)}, 24(2):751--874, 2018.

\bibitem[BW18]{BorodinWheeler}
Alexei Borodin and Michael Wheeler.
\newblock Colored stochastic vertex models and their spectral theory, 2018.
\newblock Preprint, \texttt{arXiv:1808.01866}.

\bibitem[BWo17]{Bryc}
W\l~odek Bryc and Jacek Weso\l~owski.
\newblock Asymmetric simple exclusion process with open boundaries and
  quadratic harnesses.
\newblock {\em J. Stat. Phys.}, 167(2):383--415, 2017.

\bibitem[Can16]{Cantini2}
Luigi Cantini.
\newblock Inhomogenous multispecies {TASEP} on a ring with spectral parameters,
  2016.
\newblock arXiv:1602.07921.

\bibitem[Can17]{Cantini}
L.~Cantini.
\newblock Asymmetric simple exclusion process with open boundaries and
  {K}oornwinder polynomials.
\newblock {\em Ann. Henri Poincar\'e}, 18(4):1121--1151, 2017.

\bibitem[CdGW15]{CGW}
Luigi Cantini, Jan de~Gier, and Michael Wheeler.
\newblock Matrix product formula for {M}acdonald polynomials.
\newblock {\em J. Phys. A}, 48(38):384001, 25, 2015.

\bibitem[CGdGW16]{CantiniGarbali}
Luigi Cantini, Alexandr Garbali, Jan de~Gier, and Michael Wheeler.
\newblock Koornwinder polynomials and the stationary multi-species asymmetric
  exclusion process with open boundaries.
\newblock {\em J. Phys. A}, 49(44):444002, 23, 2016.

\bibitem[CK21]{CorwinKnizel}
Ivan Corwin and Alisa Knizel.
\newblock Stationary measure for the open {KPZ} equation, 2021.
\newblock Preprint, \texttt{arXiv:2103.12253}.

\bibitem[CKS16]{CKS}
Sylvie Corteel, Jang~Soo Kim, and Dennis Stanton.
\newblock Moments of orthogonal polynomials and combinatorics.
\newblock In {\em Recent trends in combinatorics}, volume 159 of {\em IMA Vol.
  Math. Appl.}, pages 545--578. Springer, [Cham], 2016.

\bibitem[CMW22]{CMW2}
Sylvie Corteel, Olya Mandelshtam, and Lauren Williams.
\newblock From multiline queues to {M}acdonald polynomials via the exclusion
  process, 2022.
\newblock to appear in Amer. J. Math.

\bibitem[Cor12]{Corwin}
Ivan Corwin.
\newblock The {K}ardar-{P}arisi-{Z}hang equation and universality class.
\newblock {\em Random Matrices Theory Appl.}, 1(1):1130001, 76, 2012.

\bibitem[CRV14]{Crampe}
N.~Crampe, E.~Ragoucy, and M.~Vanicat.
\newblock Integrable approach to simple exclusion processes with boundaries.
  {R}eview and progress.
\newblock {\em J. Stat. Mech. Theory Exp.}, (11):P11032, 42, 2014.

\bibitem[CS14]{CS}
Luigi Cantini and Andrea Sportiello.
\newblock A one-parameter refinement of the {R}azumov-{S}troganov
  correspondence.
\newblock {\em J. Combin. Theory Ser. A}, 127:400--440, 2014.

\bibitem[CS18]{CorwinShen}
Ivan Corwin and Hao Shen.
\newblock Open {ASEP} in the weakly asymmetric regime.
\newblock {\em Comm. Pure Appl. Math.}, 71(10):2065--2128, 2018.

\bibitem[CSSW12]{CSSW}
S.~Corteel, R.~Stanley, D.~Stanton, and L.~Williams.
\newblock Formulae for {A}skey-{W}ilson moments and enumeration of staircase
  tableaux.
\newblock {\em Trans. Amer. Math. Soc.}, 364(11):6009--6037, 2012.

\bibitem[CST18]{CorwinShenTsai}
Ivan Corwin, Hao Shen, and Li-Cheng Tsai.
\newblock {${\rm ASEP}(q,j)$} converges to the {KPZ} equation.
\newblock {\em Ann. Inst. Henri Poincar\'{e} Probab. Stat.}, 54(2):995--1012,
  2018.

\bibitem[CW07a]{CW2}
Sylvie Corteel and Lauren Williams.
\newblock A {M}arkov chain on permutations which projects to the {PASEP}.
\newblock {\em Int. Math. Res. Not. IMRN}, (17):Art. ID rnm055, 27, 2007.

\bibitem[CW07b]{CW1}
Sylvie Corteel and Lauren Williams.
\newblock Tableaux combinatorics for the asymmetric exclusion process.
\newblock {\em Adv. in Appl. Math.}, 39(3):293--310, 2007.

\bibitem[CW10]{CWPNAS}
Sylvie Corteel and Lauren~K. Williams.
\newblock Staircase tableaux, the asymmetric exclusion process, and
  {A}skey-{W}ilson polynomials.
\newblock {\em Proc. Natl. Acad. Sci. USA}, 107(15):6726--6730, 2010.

\bibitem[CW11]{CW4}
Sylvie Corteel and Lauren~K. Williams.
\newblock Tableaux combinatorics for the asymmetric exclusion process and
  {A}skey-{W}ilson polynomials.
\newblock {\em Duke Math. J.}, 159(3):385--415, 2011.

\bibitem[CW18]{CW5}
Sylvie Corteel and Lauren~K. Williams.
\newblock Macdonald-{K}oornwinder moments and the two-species exclusion
  process.
\newblock {\em Selecta Math. (N.S.)}, 24(3):2275--2317, 2018.

\bibitem[DEHP93]{Derrida1}
B.~Derrida, M.~R. Evans, V.~Hakim, and V.~Pasquier.
\newblock Exact solution of a {$1$}{D} asymmetric exclusion model using a
  matrix formulation.
\newblock {\em J. Phys. A}, 26(7):1493--1517, 1993.

\bibitem[DF04]{PDF}
P.~Di~Francesco.
\newblock A refined {R}azumov-{S}troganov conjecture.
\newblock {\em J. Stat. Mech. Theory Exp.}, (8):009, 16, 2004.

\bibitem[dGR04]{DGR}
Jan de~Gier and Vladimir Rittenberg.
\newblock Refined {R}azumov-{S}troganov conjectures for open boundaries.
\newblock {\em J. Stat. Mech. Theory Exp.}, (9):009, 14, 2004.

\bibitem[DR12]{Diaconis}
Persi Diaconis and Arun Ram.
\newblock A probabilistic interpretation of the {M}acdonald polynomials.
\newblock {\em Ann. Probab.}, 40(5):1861--1896, 2012.

\bibitem[DS05]{DS}
Enrica Duchi and Gilles Schaeffer.
\newblock A combinatorial approach to jumping particles.
\newblock {\em J. Combin. Theory Ser. A}, 110(1):1--29, 2005.

\bibitem[GNS21]{Nestoridi}
Nina Gantert, Evita Nestoridi, and Dominik Schmid.
\newblock Mixing times for the simple exclusion process with open boundaries,
  2021.
\newblock Preprint, \texttt{arXiv:2003.03781}.

\bibitem[Hen72]{Hendricks}
W.~J. Hendricks.
\newblock The stationary distribution of an interesting {M}arkov chain.
\newblock {\em J. Appl. Probability}, 9:231--233, 1972.

\bibitem[Hil66]{Hill}
T.L. Hill.
\newblock Studies in irreversible thermodynamics iv. diagrammatic
  representation of steady state fluxes for uni-molercular systems.
\newblock {\em J. Theoret. Biol.}, 10:442--459, 1966.

\bibitem[ISV87]{ISV}
Mourad E.~H. Ismail, Dennis Stanton, and G\'{e}rard Viennot.
\newblock The combinatorics of {$q$}-{H}ermite polynomials and the
  {A}skey-{W}ilson integral.
\newblock {\em European J. Combin.}, 8(4):379--392, 1987.

\bibitem[Joh00]{Johansson}
Kurt Johansson.
\newblock Shape fluctuations and random matrices.
\newblock {\em Comm. Math. Phys.}, 209(2):437--476, 2000.

\bibitem[Koo92]{Koornwinder}
Tom~H. Koornwinder.
\newblock Askey-{W}ilson polynomials for root systems of type {$BC$}.
\newblock In {\em Hypergeometric functions on domains of positivity, {J}ack
  polynomials, and applications ({T}ampa, {FL}, 1991)}, volume 138 of {\em
  Contemp. Math.}, pages 189--204. Amer. Math. Soc., Providence, RI, 1992.

\bibitem[KPZ86]{KPZ}
M.~Kardar, G.~Parisi, and Y.~Zhang.
\newblock Dynamic scaling of growing interfaces.
\newblock {\em Phys. Rev.}, (56), 1986.

\bibitem[KS60]{KemenySnell}
John~G. Kemeny and J.~Laurie Snell.
\newblock {\em Finite {M}arkov chains}.
\newblock The University Series in Undergraduate Mathematics. D. Van Nostrand
  Co., Inc., Princeton, N.J.-Toronto-London-New York, 1960.

\bibitem[KW21]{KW}
Donghyun Kim and Lauren Williams.
\newblock Schubert polynomials, the inhomogeneous {TASEP}, and evil-avoiding
  permutations, 2021.
\newblock arXiv:2106.13378.

\bibitem[Lig85]{Liggett}
Thomas~M. Liggett.
\newblock {\em Interacting particle systems}, volume 276 of {\em Grundlehren
  der mathematischen Wissenschaften [Fundamental Principles of Mathematical
  Sciences]}.
\newblock Springer-Verlag, New York, 1985.

\bibitem[LP16]{Lyons}
Russell Lyons and Yuval Peres.
\newblock {\em Probability on trees and networks}, volume~42 of {\em Cambridge
  Series in Statistical and Probabilistic Mathematics}.
\newblock Cambridge University Press, New York, 2016.

\bibitem[LR83]{Leighton}
Frank~Thomson Leighton and Ronald~L. Rivest.
\newblock Estimating a probability using finite memory.
\newblock In {\em Foundations of computation theory ({B}orgholm, 1983)}, volume
  158 of {\em Lecture Notes in Comput. Sci.}, pages 255--269. Springer, Berlin,
  1983.

\bibitem[LW12]{LW}
Thomas Lam and Lauren Williams.
\newblock A {M}arkov chain on the symmetric group that is {S}chubert positive?
\newblock {\em Exp. Math.}, 21(2):189--192, 2012.

\bibitem[Mac95]{Macdonald}
I.~G. Macdonald.
\newblock {\em Symmetric functions and {H}all polynomials}.
\newblock Oxford Mathematical Monographs. The Clarendon Press, Oxford
  University Press, New York, second edition, 1995.

\bibitem[Mar20]{Martin}
James~B. Martin.
\newblock Stationary distributions of the multi-type {ASEP}.
\newblock {\em Electron. J. Probab.}, 25:Paper No. 43, 41, 2020.

\bibitem[MGP68]{bio}
J.~Macdonald, J.~Gibbs, and A.~Pipkin.
\newblock Kinetics of biopolymerization on nucleic acid templates.
\newblock {\em Biopolymers}, 6, 1968.

\bibitem[Pan19]{Pang}
C.~Y.~Amy Pang.
\newblock Lumpings of algebraic {M}arkov chains arise from subquotients.
\newblock {\em J. Theoret. Probab.}, 32(4):1804--1844, 2019.

\bibitem[PT18]{Pitman}
Jim Pitman and Wenpin Tang.
\newblock Tree formulas, mean first passage times and {K}emeny's constant of a
  {M}arkov chain.
\newblock {\em Bernoulli}, 24(3):1942--1972, 2018.

\bibitem[Qua12]{Quastel}
Jeremy Quastel.
\newblock Introduction to {KPZ}.
\newblock In {\em Current developments in mathematics, 2011}, pages 125--194.
  Int. Press, Somerville, MA, 2012.

\bibitem[San94]{Sandow}
S.~Sandow.
\newblock Partially asymmetric exclusion process with open boundaries.
\newblock {\em Phys. Rev.}, (E50), 1994.

\bibitem[Spi70]{Spitzer}
Frank Spitzer.
\newblock Interaction of {M}arkov processes.
\newblock {\em Advances in Math.}, 5:246--290, 1970.

\bibitem[Sta99]{EC2}
Richard~P. Stanley.
\newblock {\em Enumerative combinatorics. {V}ol. 2}, volume~62 of {\em
  Cambridge Studies in Advanced Mathematics}.
\newblock Cambridge University Press, Cambridge, 1999.
\newblock With a foreword by Gian-Carlo Rota and appendix 1 by Sergey Fomin.

\bibitem[Tse63]{Tsetlin}
M.L. Tsetlin.
\newblock Finite automata and models of simple forms of behavior.
\newblock {\em Russian Mathematical Surveys}, 18:1--28, 1963.

\bibitem[TW09]{TracyWidom}
Craig~A. Tracy and Harold Widom.
\newblock Asymptotics in {ASEP} with step initial condition.
\newblock {\em Comm. Math. Phys.}, 290(1):129--154, 2009.

\bibitem[USW04]{USW}
Masaru Uchiyama, Tomohiro Sasamoto, and Miki Wadati.
\newblock Asymmetric simple exclusion process with open boundaries and
  {A}skey-{W}ilson polynomials.
\newblock {\em J. Phys. A}, 37(18):4985--5002, 2004.

\bibitem[Vie85]{Viennot-book}
G\'{e}rard Viennot.
\newblock A combinatorial theory for general orthogonal polynomials with
  extensions and applications.
\newblock In {\em Orthogonal polynomials and applications ({B}ar-le-{D}uc,
  1984)}, volume 1171 of {\em Lecture Notes in Math.}, pages 139--157.
  Springer, Berlin, 1985.

\end{thebibliography}

\end{document}